\documentclass[10pt]{amsart}
\usepackage[dvips]{graphicx}

\newtheorem{theorem}{Theorem}[section]
\newtheorem{lemma}[theorem]{Lemma}

\newtheorem{prop}[theorem]{Proposition}

\theoremstyle{definition}

\theoremstyle{remark}

\numberwithin{equation}{section}

\title[Discrete Schr{\"o}dinger Operators]{Inverse problems,
trace formulae for Discrete Schr{\"o}dinger Operators}
\author{Hiroshi ISOZAKI}
\address{Institute of Mathematics,
University of Tsukuba
\\
Tsukuba, 305-8571, JAPAN
\ isozakih@math.tsukuba.ac.jp}
\author{Evgeny Korotyaev}
\address{Leningrad State
University, Peterusburgskoe shosse 10,
Saint-Petersburg, 196605, Russia
 \ korotyaev@gmail.com}

\date{March, 7, 2011}

\begin{document}

\maketitle

\begin{abstract}
We study discrete Schr{\"o}dinger operators with compactly
supported potentials on ${\bf Z}^d$. Constructing spectral
representations and representing S-matrices by the generalized
eigenfunctions, we show that the potential is uniquely reconstructed
from the S-matrix of all energies. We also study the spectral shift
function $\xi(\lambda)$ for the trace class potentials, and estimate
the discrete spectrum in terms of the moments of $\xi(\lambda)$ and
the potential.
\end{abstract}


\section{Introduction}

In this paper, we consider trace formulas and inverse scattering problems
for Schr{\"o}dinger operators on the square lattice ${\bf Z}^d$ with
$d \geq 2$. We restrict ourselves to the case of compactly supported,
or trace class, potentials. Our aim is two-fold :  the recontsruction
of the potential from the scattering matrix, and the computation of
trace formula using the spectral shift function. We begin with
the forward problem. We shall construct the generalized Fourier transform
and represent the S-matrix by generalized eigenfunctions.
We then show that given the S-matrix $\mathcal S(\lambda)$ for all
energies, one can uniquely  reconstruct the potential (Theorem 4.4).
We next compute the asymptotic expansion of the perturbation
determinant associated with the discrete Hamiltonian.
By virtue of Krein's spectral shift function, one can compute
the moments of $\log \det \mathcal S(\lambda)$. As a by-product,
one can estimate the discrete spectrum using these moments (Theorem 6.4).

In the continuous model, the first mathematical result on the inverse
scattering for multi-dimensional Schr{\"o}dinger operators was that of
Faddeev \cite{Fa56} : the reconstruction of the potential from
the high-energy behavior of the scattering matrix using
the Born approximation. In this paper, we give its discrete analogue.
Instead of high-energy, we consider the analytic continuation of
the S-matrix with respect to the energy parameter and use
the complex Born approximation. The analytic property of
the resolvent of the discrete Laplacian is more complicated than
the continuous case, which requires harder analysis in studying
the inverse scattering problem. In the continuous case,
Faddeev proposed a multi-dimensional analogue of
the Gel'fand-Levitan theory using new Green's function of
the Helmholtz equation (\cite{Fa66}, \cite{Fa76}).
Faddeev's Green function was rediscovered and developed in 1980's
by Sylvester-Uhlmann \cite{SyUh87}, Nachman \cite{Na88},
Khenkin-Novikov \cite{KhNo87} (see the survey article \cite{Is03}).
Reconstruction of the potential from the S-matrix of a fixed energy
is one of the novelties of this idea. In the forthcoming paper,
we shall study the inverse scattering from one fixed energy in
the discrete model.

Our next purpose is the trace formula.
 It is well-known that the scattering matrix, Krein's spectral shift
 function and the perturbation determinant for a pair of self-adjoint
 operators $\widehat H_0$ and $\widehat H = \widehat H_0 + \widehat V$
 are mutually related. We shall write down the first 5 moments of
 the spectral shift function for our discrete model in terms of
 the traces of the potential $\widehat V$. If the potential admits
 a definite sign, we can obtain estimates of the discrete spectrum
 by the traces of $\widehat V$.

The computation of trace constitutes the central part of the study
of spectral theory, since it provides quantitative information of
the operator in question, hence serves as a clue to the inverse problem.
In the continuous case, the trace formula was first obtained by
Buslaev-Faddeev \cite{BF} in the one-dimensional problem and by
Buslaev \cite{B1}, \cite{B2} in the three-dimensional problem. Since then,
an abundance of articles have been devoted to this subject, see e.g.
\cite{CV} and \cite{G}, \cite{P}, \cite{R}, \cite{GHS95}.
Gesztesy-Holden-Simon-Zhao \cite{GHSZ96} computed the trace
${\rm Tr}\,(L - L_A)$, where $L = - \Delta + V$ is a Schr{\"o}dinger
operator for the continuous case, and $L_A$ is $L$ with Dirichlet
condition on a subset $A \subset {\bf R}^{\nu}$. Shirai \cite{Sh}
studied this problem on a graph.
 Karachalios \cite{Ka} and Rosenblum-Solomjak  \cite{RoSol09} computed
 the Cwikel-Lieb-Rosenblum type bound for the discrete Schr{\"o}dinger
 operator. The well-known Effimov effect has a different property in
 the case of the discrete model. See e.g. Albeverio, Dell Antonio and
 Lakaev \cite{AAL07}. See also \cite{ALMM06}.

In \S 2, we shall prove the limiting absorption principle with the aid of
Mourre's commutator theory \cite{Mo81}. We then derive the spectral
representation in \S 3, and represent the S-matrix by generalized
eigenfunctions. The inverse scattering problem is solved in \S 4,
and \S 6 is devoted to the trace formula.

The essential spectrum of our discrete Schr{\"o}dinger operator
$\widehat H = \widehat H_0 +\widehat V$ on ${\bf Z}^d$ fills
the interval $[0,d]$, and the set $[0,d]\cap{\bf Z}$ is that
of critical points, since $\widehat H_0$ is unitarily equivalent
to the operator of multiplication by $(d - \sum_{j=1}^d\cos x_j)/2$ on
the torus ${\bf R}^d/(2\pi{\bf Z})^d$. In fact, the resolvent estimates
in \S 2 are proved outside $[0,d]\cap{\bf Z}$. The behavior of the free
resolvent $(\widehat H_0-z)^{-1}$ near the critical points depends on
the dimension $d$. This is itself interesting and is studied in \S 5
(Lemmas 5.3, 5.4 and 5.5), although we do not use it in this paper.

The notation used in this paper is standard. For two Banach spaces $X$
and $Y$, ${\bf B}(X;Y)$ denotes the set of all bounded operators from
$X$ to $Y$. For a self-adjoint operator $A$, $\sigma(A)$, $\sigma_{p}(A)$,
$\sigma_{d}(A)$, $\sigma_{e}(A)$, $\sigma_{ac}(A)$ and $\rho(A)$ denote its
spectrum, point spectrum (= the set of all eigenvalues), discrete spectrum,
 essential spectrum, absolutely continuous spectrum and resolvent set,
 respectively. For a trace class operator $K$, ${\rm Tr}\,(K)$ denotes
 the trace of $K$.


\section{Schr{\"o}dinger operators on the lattice}


\subsection{Discrete Schr{\"o}dinger operator}
Let ${\bf Z}^d = \{n = (n_1,\cdots,n_d) \, ; \, n_i \in {\bf Z}\}$,
and $e_{1} = (1,0,\cdots,0), \cdots, e_{d} = (0,\cdots,0,1)$ be the standard
bases of ${\bf Z}^d$.  The  Schr{\"o}dinger operator
${\widehat H}$ on ${\bf Z}^d$ is defined by
\begin{equation}
{\widehat H} = {\widehat H}_{0} + {\widehat V},
\nonumber
\end{equation}
where for ${\widehat f} = \big\{\hat f(n)\big\}_{n\in{\bf Z}^d}
\in l^{2}({\bf Z}^d)$ and $n \in {\bf Z}^d$
\begin{equation}
\big({\widehat H}_{0}{\widehat f}\big)(n) = \frac{d}{2}\,{\widehat f}(n) -
\frac{1}{4}\sum_{j=1}^{d}\big\{{\widehat f}(n + e_{j}) +
{\widehat f}(n - e_{j})\big\},
\nonumber
\end{equation}
\begin{equation}
({\widehat V}{\widehat f})(n) = \widehat V(n){\widehat f}(n).
\nonumber
\end{equation}
Until the end of \S 4, we impose the following assumption on $\widehat V$:

\bigskip
\noindent
{\bf (A-1)}  $\ \widehat V$ {\it is real-valued, and
$\widehat V(n) = 0$ except for a finite number of} $n$.

\bigskip

 Define the 1-dimensional projection $\widehat P(n)$ by
\begin{equation}
\left(\widehat P(n)\widehat f\right)(m) = \delta_{nm}\widehat f(m).
\nonumber
\end{equation}
Then $\widehat V$ is rewritten as
\begin{equation}
\widehat V = \sum_n\widehat V(n)\widehat P(n).
\nonumber
\end{equation}
Let us introduce the shift operator
\begin{equation}
\big(\widehat S_j\widehat f\big)(n) = \widehat f(n + e_j),
\quad
\big((\widehat S_j)^{\ast}\widehat f\big)(n) = \widehat f(n - e_j),
\nonumber
\end{equation}
and the position operator
\begin{equation}
\big(\widehat N_j\widehat f\big)(n) = n_j\widehat f(n).
\nonumber
\end{equation}
A direct computation yields the following lemma.


\begin{lemma} $\ \widehat S_j$ is unitary, $\widehat N_j$ is self-adjoint
with its natural domain and
$$
\big[\widehat S_j,\widehat N_j\big] = \widehat S_j.
$$
\end{lemma}

Letting $\widehat N = (\widehat N_1,\cdots,\widehat N_d)$, $\widehat H$
is rewritten as
\begin{equation}
\widehat H =  \frac{d}{2} - \frac{1}{4}\sum_{j=1}^d\left(\widehat S_j +
(\widehat S_j)^{\ast}\right)  + \widehat V(\widehat N).
\nonumber
\end{equation}

 The spectral properties of ${\widehat H}$ is easier to  describe by passing
to its unitary transformation by the Fourier series. Let
\begin{equation}
{\bf T}^d = {\bf R}^{d}/(2\pi{\bf Z})^{d} = [-\pi,\pi]^d
\nonumber
\end{equation}
be the flat torus and ${\mathcal U}$ the unitary operator from
$l^{2}({\bf Z}^d)$ to $L^{2}({\bf T}^d)$ defined by
\begin{equation}
({\mathcal U}\,{\widehat f})(x) = (2\pi)^{-d/2}
\sum_{n\in{\bf Z}^d}{\widehat f}(n)
e^{-in\cdot x}.
\nonumber
\end{equation}
The shift operator and the position operator are rewritten as
\begin{equation}
S_j := \mathcal U\,\widehat S_j\,\mathcal U^{\ast} = e^{ix_j}, \quad
N_j := \mathcal U\, \widehat N_j\, \mathcal U^{\ast} = i\partial/\partial x_j.
\nonumber
\end{equation}
Letting
\begin{equation}
H_0 = {\mathcal U}\, {\widehat H_0}\, {\mathcal U}^{\ast}, \quad
V = {\mathcal U}\, {\widehat V}\, {\mathcal U}^{\ast},
\nonumber
\end{equation}
we have
\begin{equation}
H := {\mathcal U}{\widehat H}{\mathcal U}^{\ast} = H_0 + V,
\label{S2DefH}
\end{equation}
\begin{equation}
H_0 = \frac{1}{2}\Big(d - \sum_{j=1}^d\cos x_j\Big) =: h(x),
\label{S2Defp(x)}
\end{equation}
\begin{equation}
(Vf)(x) = (2\pi)^{-d/2}\int_{{\bf T}^d}V(x - y)f(y)dy,
\nonumber
\end{equation}
\begin{equation}
V(x) = (2\pi)^{-d/2}\sum_{n\in{\bf Z}^d}{\widehat V}(n)e^{-in\cdot x}.
\nonumber
\end{equation}
In fact, this is a special case of the Friedrichs model (see e.g. \cite{Fa64}).
The following theorem is easily proven by (\ref{S2DefH}) and Weyl's theorem.


\begin{theorem}
(1) $\  \sigma(H_{0}) = \sigma_{ac}(H_{0}) = [0,d]$. \\
\noindent
(2) $\  $ $\sigma_{e}(H) = [0,d], \quad
\sigma_{d}(H) \subset {\bf R}\setminus[0,d]$.
\end{theorem}


\subsection{Sobolev and Besov spaces}
We put $N = (N_1,\cdots,N_d)$, and let $N^2$ be the self-adjont operator
defined by
\begin{equation}
N^2 = \sum_{j=1}^dN_j^2 =  - \Delta, \quad {\rm on} \quad {\bf T}^d,
\nonumber
\end{equation}
where $\Delta$ denotes the Laplacian on ${\bf T}^d = [-\pi,\pi]^d$ with
periodic boundary condition.  We put
\begin{equation}
|N| = \sqrt{N^2} = \sqrt{-\Delta}.
\nonumber
\end{equation}
We introduce the norm
\begin{equation}
\|u\|_s = \|(1 + N^2)^{s/2}u\|, \quad s \in {\bf R},
\nonumber
\end{equation}
$\|\cdot\|$ being the norm on $L^2({\bf T}^d)$,
and let ${\mathcal H}^s$ be the completion of $D(|N|^s)$, the domain of
$|N|^s$, with respect to the norm
$\|u\|_{s}$ :
\begin{equation}
{\mathcal H}^s = \{u \in {\mathcal D}^{\prime}({\bf T}^d)\, ; \,
\|u\|_{s} < \infty\}, \nonumber
\end{equation}
where $\mathcal D^{\prime}({\bf T}^d)$ denotes the space of distribution
on ${\bf T}^d$.
Put $\mathcal H = \mathcal H^0 = L^2({\bf T}^d)$.

For a self-adjoint operator $T$, let $\chi(a \leq T < b)$ denote the operator
$\chi_{I}(T)$, where $\chi_I(\lambda)$ is the characteristic function of
the interval $I = [a,b)$. The operators $\chi(T < a)$ and $\chi(T \geq b)$
are defined similarly.  Using the sequence $\{r_j\}_{j=0}^{\infty}$ with
$r_{-1} = 0$, $r_j = 2^j \ (j \geq 0)$, we define the Besov space
$\mathcal B$ by
\begin{equation}
\mathcal B = \Big\{f \in {\mathcal H}\, ; \|f\|_{\mathcal B}
=\sum_{j=0}^{\infty}r_j^{1/2}\|\chi(r_{j-1} \leq |N| < r_j)f\| <
\infty\Big\}. \nonumber
\end{equation}
Its dual space $\mathcal B^{\ast}$ is the completion of $\mathcal H$ by
the following norm
\begin{equation}
\|u\|_{\mathcal B^{\ast}}= \sup_{j\geq 0}r_j^{-1/2}\|\chi(r_{j-1}
\leq |N|< r_j)u\|.\nonumber
\end{equation}

The following Lemmas 2.3 and 2.4 are proved in the same way as in
the continuous case
(\cite{AgHo76}). We omit the proof.


\begin{lemma}
There exists a constant $C > 0$ such that
\begin{equation}
C^{-1}\|u\|_{\mathcal B^{\ast}}\leq \left(\sup_{R>1}\frac{1}{R}\|\chi(|N|
< R)u\|^2\right)^{1/2}\leq C\|u\|_{\mathcal B^{\ast}}.
\nonumber
\end{equation}
\end{lemma}

Therefore, in the following, we use
\begin{equation}
\|u\|_{\mathcal B^{\ast}} = \left(\sup_{R>1}\frac{1}{R}\|\chi(|N|
<R)u\|^2\right)^{1/2} \nonumber
\end{equation}
as the norm on $\mathcal B^{\ast}$.


\begin{lemma}
For $s > 1/2$, the following inclusion relations hold :
\begin{equation}
\mathcal H^{s} \subset \mathcal B \subset \mathcal H^{1/2} \subset
\mathcal H \subset \mathcal H^{-1/2} \subset \mathcal B^{\ast} \subset
\mathcal H^{-s}.
\nonumber
\end{equation}
\end{lemma}

\medskip
We also put $\widehat{\mathcal H} = l^2({\bf Z}^d)$, and define
$\widehat{\mathcal H}^s$, $\widehat{\mathcal B}$,
$\widehat{\mathcal B}^{\ast}$ by replacing $N$ by $\widehat N$.
Note that $\widehat{\mathcal H}^s = \mathcal U^{\ast}\mathcal H^s$
and so on. In particular,
Parseval's formula implies that
\begin{equation}
\|u\|_{\mathcal H^s}^2 = \|\widehat u\|_{\widehat{\mathcal H}^s}^2 =
 \sum_{n\in Z^d}(1 + |n|^2)^s|\widehat u(n)|^2 ,
\nonumber
\end{equation}
\begin{equation}
\|u\|_{{\mathcal B}^{\ast}}^2 =
\|\widehat u\|_{\widehat{\mathcal B}^{\ast}}^2
=
\sup_{R>1}\frac{1}{R}\sum_{|n|<R}|\widehat u(n)|^2,
\nonumber
\end{equation}
$\widehat u(n)$ being the Fourier coefficient of $u(x)$.


\subsection{Mourre estimate}
 Let
$\widehat{\mathcal H}^{\infty} = \cap_{s>0}\widehat{\mathcal H}^s$,
and  define a symmetric operator $\widehat A$ with domain
$\widehat{\mathcal H}^{\infty}$ by
\begin{equation}
\widehat A = i \sum_{j=1}^d\left(\big((\widehat S_j)^{\ast} -
\widehat S_j\big)\widehat N_j + \frac{\widehat S_j +
(\widehat S_j)^{\ast}}{2}\right).
\label{S2WidehatA}
\end{equation}
Then $\widehat A$ is essentially self-adjoint. In fact, letting
$\widehat M = 1 + \widehat N^2$, we can find a constant $C > 0$ such that
\begin{equation}
\|\widehat A\widehat u\| \leq C\|\widehat M\widehat u\|,
\quad
|(\widehat A\widehat u, \widehat M\widehat u) - (\widehat M\widehat u,
\widehat A\widehat u)| \leq C\|\widehat M^{1/2}\widehat u\|^2, \quad
\forall \widehat u \in \widehat{\mathcal H}^{\infty}.
\end{equation}
 Nelson's commutator theorem (\cite{ReSi75}, p. 193) then implies the result.

By Lemma 2.1, (\ref{S2WidehatA}) is rewritten as
\begin{equation}
\widehat A = -i\sum_{j=1}^d\left(\widehat N_j\widehat S_j -
(\widehat S_j)^{\ast}\widehat N_j  + \frac{\widehat S_j -
(\widehat S_j)^{\ast}}{2}\right).
\nonumber
\end{equation}
Let us note that in \cite{AMBJS99}, $- i\sum_{j=1}^d
\left(\widehat N_j\widehat S_j - (\widehat S_j)^{\ast}
\widehat N_j\right)$ is used as $\widehat A$.
Our choice of $\widehat A$ comes from the following reasoning. Let
$h(x)$ be defined by (\ref{S2Defp(x)}).
Passing to the Fourier series, we have
\begin{equation}
A = \mathcal U\, \widehat A\, \mathcal U^{\ast} = i\sum_{j=1}^d
\Big(2\sin x_j\frac{\partial}{\partial x_j} + \cos x_j\Big) =
2i\big(\nabla_xh\cdot\nabla_x + \nabla_x\cdot(\nabla_xh)\big).
\nonumber
\end{equation}
This is an analogue of the generator of dilation group on ${\bf R}^d$.
 We then have
\begin{equation}
i[H_{0},A] = 4|\nabla_x h|^2 = \sum_{j=1}^d(\sin x_j)^2.
\nonumber
\end{equation}
Let $E_{H_0}(\lambda)$ and $E_H(\lambda)$ be the spectral decompositions
of $H_0$ and $H$, respectively.


\begin{lemma} Let $\lambda \in (0,d)\setminus {\bf Z}$. Then there exist
constants $\delta, C > 0$ and a compact operator $K$ such that
\begin{equation}
E_{H}(I)[H,iA]E_{H}(I) \geq CE_H(I) + K, \quad I =
(\lambda-\delta,\lambda+\delta).
\nonumber
\end{equation}
\end{lemma}

Proof.
For $\lambda \in (0,d)\setminus {\bf Z}$, let
\begin{equation}
M_{\lambda} = \{x \in {\bf T}^d\, ; \, h(x) = \lambda\}.
\nonumber
\end{equation}
If $\nabla h(x) = 0$, then $\cos x_j = \pm 1$, and $h(x) \in {\bf Z}$.
Therefore, the assumption $\lambda \not\in {\bf Z}$ implies that
$\nabla h(x) \neq 0$ on $M_{\lambda}$, hence $M_{\lambda}$ is a
real analytic manifold.  We put
\begin{equation}
C_{0}(\lambda) = \inf_{x \in M_{\lambda}}|\nabla h(x)|^2.
\nonumber
\end{equation}
Then for any small $\epsilon > 0$, there exists $\delta > 0$ such that
\begin{equation}
|\nabla h(x)|^{2} \geq C_{0}(\lambda) - \epsilon \quad {\rm on} \quad
h^{-1}([\lambda - \delta,\lambda + \delta]).
\nonumber
\end{equation}
We have, therefore,
\begin{equation}
E_{H_{0}}(I)[H_{0},iA]E_{H_{0}}(I) \geq (C_{0}(\lambda) - \epsilon)E_{H_{0}}(I),
\quad I = (\lambda-\delta,\lambda+\delta).
\nonumber
\end{equation}
Since $\widehat V$ is a compact operator, so are $V$ and $E_{H}(I) - E_{H_0}(I)$. This proves the lemma. \qed

\medskip
Let $R(z) = (H - z)^{-1}$ be the resolvent of $H$.


\begin{theorem}
(1) $\sigma_{p}(H)\cap \big((0,d)\setminus {\bf Z}\big)$ is discrete and
finite multiplicities with possible accumulation points in ${\bf Z}$. \\
\noindent
(2) Let $s > 1/2$ and $\lambda \in (0,d)\setminus\big({\bf Z}\cup
\sigma_p(H)\big)$.  Then, thers exists a norm limit $R(\lambda \pm i0) :=
 \lim_{\epsilon\to0} R(\lambda \pm i\epsilon) \in {\bf B}
 (\mathcal H^s;\mathcal H^{-s})$. Moreover, we have
\begin{equation}
\sup_{\lambda \in J} \|R(\lambda \pm i0)\|_{{\bf B}
(\mathcal B;\mathcal B^{\ast})} < \infty,
\label{S2LAP}
\end{equation}
for any  compact interval $J$ in $(0,d)\setminus\big
({\bf Z}\cup\sigma_p(H)\big)$. The mapping $(0,d)
\setminus\big({\bf Z}\cup\sigma_p(H)\big) \ni \lambda \to
R(\lambda \pm i0)$ is norm continuous in ${\bf B}({\mathcal H^s};
\mathcal H^{-s})$ and weakly continuous in ${\bf B}
(\mathcal B\,;\mathcal B^{\ast})$.\\
\noindent
(3) $H$ has no singular continuous spectrum.
\end{theorem}

This theorem follows from the well-known Mourre theory. We shall
give here a brief explanation.
First we introduce an abstract Besov space. We define
\begin{equation}
\mathcal B_A = \Big\{f \in \mathcal H\, ; \, \|f\|_{\mathcal B_A} =
\sum_{j=0}^{\infty}r_j^{1/2}\|\chi(r_{j-1}\leq |A| < r_j)f\| < \infty\Big\},
\nonumber
\end{equation}
where $\mathcal H = L^2({\bf T}^d)$. Its dual space
${\mathcal B_A}^{\ast}$ is the completion of $\mathcal H$ by the norm
\begin{equation}
\|u\|_{{\mathcal B_A}^{\ast}} = \sup_j r_j^{-1/2}\|\chi(r_{j-1}
\leq |A| < r_j)u\|.
\nonumber
\end{equation}
The abstract theory of Mourre based on Lemma 2.5 then yields


\begin{lemma} Let $J$ be as in Theorem 2.6 (2). Then there exists a
constant $C > 0$ such that
\begin{equation}
\sup_{{\rm Re}\, z \in J, {\rm Im}\, z
\neq 0}\|(H - z)^{-1}f\|_{{\mathcal B_A}^{\ast}} \leq C\|f\|_{\mathcal B_A},
\quad
\forall f \in \mathcal B_A.
\nonumber
\end{equation}
\end{lemma}

For the proof of the lemma, see \cite{ABG96}, \cite{AMBJS99}, \cite{JePe85}.
Therefore, to prove Theorem 2.6, we have only to replace $\mathcal B_A$ by
 $\mathcal B$ using the following lemma.


\begin{lemma}
There is a constant $C > 0$ such that
\begin{equation}
\|f\|_{\mathcal B_A} \leq C\|f\|_{\mathcal B}, \quad \forall f \in \mathcal B.
\nonumber
\end{equation}
\end{lemma}

Proof. For $t \in {\bf R}$, let $\langle t\rangle = (1 + t^2)^{1/2}$. By
definitions of $A$ and $N$ we have
\begin{equation}
\langle A\rangle \langle N\rangle^{-1} \in {\bf B}(\mathcal H\, ;  \mathcal H).
\nonumber
\end{equation}
For $f \in \mathcal B$, we put $f_j = \chi(r_{j-1}\leq |N| < r_j)f$. Then
\begin{equation}
\|\langle A\rangle f_j\| \leq C\|\langle N\rangle f_j\| \leq Cr_j\|f_j\|,
\nonumber
\end{equation}
which implies
\begin{equation}
\begin{split}
\|\chi(r_{k-1}\leq |A| < r_k)f_j\| &= \|\chi(r_{k-1}\leq |A| < r_k)
\langle A\rangle^{-1}\langle A\rangle f_j\| \\
& \leq Cr_k^{-1}\|\langle A \rangle f_j\| \leq Cr_k^{-1}r_j\|f_j\|.
\end{split}
\nonumber
\end{equation}
Then we have
\begin{equation}
\sum_{k>j}r_k^{1/2}\|\chi(r_{k-1}\leq |A|<r_k)f_j\| \leq
C\sum_{k>j}r_k^{-1/2}r_j\|f_j\| \leq
Cr_j^{1/2}\|f_j\|,
\nonumber
\end{equation}
\begin{equation}
\sum_{k\leq j}r_k^{1/2}\|\chi(r_{k-1}\leq |A|<r_k)f_j\| \leq
\sum_{k\leq j}r_k^{1/2}\|f_j\| \leq
Cr_j^{1/2}\|f_j\|.
\nonumber
\end{equation}
We have, therefore,
\begin{equation}
\|f_j\|_{\mathcal B_A} \leq Cr_j^{1/2}\|f_j\|, \quad j = 0,1,2,\cdots.
\nonumber
\end{equation}
Summing up these inequalities with respect to $j$, we obtain the lemma. \qed


\section{Spectral representations and S-matrices}


\subsection{Spectral representation on the torus}

For $t \in (0,d)\setminus{\bf Z}$, let $dM_t$ be the mesure on $M_t$ induced
from $dx$. By taking $t = h(x)$ as a new variable, one can show that for
$f \in C({\bf T}^d)$ supported in $\{x \in {\bf T}^d\, ; \, h(x)\not\in
{\bf Z}\}$
\begin{equation}
\int_{{\bf T}^d}f(x)dx = \int_0^d\left(\int_{M_t}f\frac{dM_t}{|\nabla_x h|}\right)dt.
\label{S3dxdMtdt}
\end{equation}
For $f, g \in L^2({\bf T}^d)$, we have
$$
\left(R_0(z)f,g\right) = \int_{{\bf T}^d}\frac{f(x)\overline{g(x)}}{h(x) - z}dx.
$$
Therefore, if $f, g \in C^1({\bf T}^d)$ and $\lambda \in (0,d)\setminus{\bf Z}$,
\begin{equation}
\left(R_0(\lambda \pm i0)f,g\right) = \pm i\pi\int_{M_{\lambda}}f\overline{g}\,
\frac{dM_{\lambda}}{|\nabla_xh |} + {\rm p.v.}\int_{{\bf T}^d}
\frac{f\overline{g}}{h(x)-\lambda}dx.
\label{S3Privalov}
\end{equation}
 Let $L^2(M_{\lambda})$ be the Hilbert space equipped with the inner product
\begin{equation}
(\varphi,\psi)_{L^2(M_{\lambda})} = \int_{M_{\lambda}}
\varphi\,\overline{\psi}\,\frac{dM_{\lambda}}{|\nabla_xh|}.
\label{S3L2Mlambdainnerprod}
\end{equation}
We define
\begin{equation}
\mathcal F_0(\lambda)f = f\Big|_{M_{\lambda}},
\label{S3trace}
\end{equation}
 where the right-hand means the trace on , i.e. the restriction to,
 $M_{\lambda}$.
Then we have by (\ref{S3Privalov})


\begin{lemma}
 For $\lambda \in (0,d)\setminus {\bf Z}$, and $f, g \in C^1({\bf T}^d)$,
\begin{equation}
 \frac{1}{2\pi i}\big((R_0(\lambda + i0) - R_0(\lambda - i0))f,g\big)
 = \big(\mathcal F_0(\lambda)f,\mathcal F_0(\lambda)g\big)_{L^2(M_{\lambda})}.
\nonumber
\end{equation}
\end{lemma}

By (\ref{S2LAP}), this lemma implies
\begin{equation}
\mathcal F_0(\lambda) \in {\bf B}(\mathcal B;L^2(M_{\lambda})).
\label{S2F0lambdabounded}
\end{equation}
Moreover,
\begin{equation}
(f,g)_{L^2(T^{d})} = \int_0^{d}(\mathcal F_0(\lambda)f,\mathcal
F_0(\lambda)g)_{L^2(M_{\lambda})}d\lambda, \quad f, g \in \mathcal B.
\label{S2F0isometry}
\end{equation}
 The adjoint operator $\mathcal F_0(\lambda)^{\ast}$ is defined by
$$
(\mathcal F_0(\lambda)f,\phi)_{L^2(M_{\lambda})} =
(f,\mathcal F_0(\lambda)^{\ast}\phi)_{L^2({\bf T}^d)}.
$$
By (\ref{S2F0lambdabounded}), $\mathcal F_0(\lambda)^{\ast}
\in {\bf B}( L^2(M_{\lambda});\mathcal B^{\ast})$, and by
(\ref{S3trace}), $\mathcal F_0(\lambda)(H_0 - \lambda) = 0$. Hence we have
\begin{equation}
(H_0 - \lambda)\mathcal F_0(\lambda)^{\ast} = 0.
\nonumber
\end{equation}

In view of (\ref{S3dxdMtdt}), we can identify $L^2({\bf T}^d)$ with
the space of $L^2$-functions $f(\lambda)$ over $(0,d)$ with respect
to the measure $d\lambda$ such that for a.e. $\lambda \in (0,d)$,
$f(\lambda)$ takes values in $L^2(M_{\lambda})$. We denote this space
by $L^2((0,d);L^2(M_{\lambda});d\lambda)$.

We put $(\mathcal F_0f)(\lambda) = \mathcal F_0(\lambda)f$ for
$f \in \mathcal B$.
The following Theorem 3.2 is essentially a reinterpretation of
the identification $L^2({\bf T}^d) \simeq L^2((0,d);L^2(M_{\lambda});
d\lambda)$. However, we give a functional analytic proof for the later
 convenience.


\begin{theorem}
(1) $\mathcal F_0$ is uniquely extended to a unitary operator
 $$
 \mathcal F_0 : L^2(T^{d}) \to L^2((0,d);L^2(M_{\lambda});d\lambda).
 $$
(2) $\mathcal F_0$ diagonalizes $H_0$ :
$$
(\mathcal F_0H_0f)(\lambda) = \lambda (\mathcal F_0f)(\lambda), \quad
\forall f \in L^2({\bf T}^d).
$$
(3) For any compact interval $I \subset (0,d)\setminus{\bf Z}$,
$$
\int_I\mathcal F_0(\lambda)^{\ast}g(\lambda)d\lambda \in L^2(T^d), \quad
\forall g \in L^2((0,d);L^2(M_{\lambda});d\lambda).
$$
Moreover, for $I_N = \cup_{j=1}^d(j-1+1/N,j-1/N)$,  the inversion formula holds:
$$
f = \lim_{N\to\infty}\int_{I_N}\mathcal F_0(\lambda)^{\ast}
(\mathcal F_0f)(\lambda)d\lambda, \quad \forall f \in L^2({\bf T}^d),
$$
where the limit is taken in the norm of $L^2({\bf T}^d)$.
\end{theorem}

Proof. By (\ref{S2F0isometry}), $\mathcal F_0$ is uniquely extended to
an isometric operator from $L^2(T^{d})$ to $ L^2((0,d);L^2(M_{\lambda});
d\lambda)$. To show that it is onto, we have only to note that the range
of $\mathcal F_0$ is dense. For $f \in \mathcal B$, we have $\mathcal
F_0(\lambda)(H_0-\lambda)f = 0$ by definition, which proves (2). To show
(3), we first note that for a compact interval $I \subset (0,d)
\setminus{\bf Z}$,
$\int_I\mathcal F_0(\lambda)^{\ast}g(\lambda)d\lambda \in \mathcal
B^{\ast}$.
We use $(\, ,\, )$ to denote the inner product of $L^2({\bf T}^d)$
as well as the coupling of $\mathcal B$ and $\mathcal B^{\ast}$. Then we have
\begin{equation}
\begin{split}
\left(\int_I\mathcal F_0(\lambda)^{\ast}g(\lambda)d\lambda,f\right) & =
\int_I\left(\mathcal F_0(\lambda)^{\ast}g(\lambda),f\right)d\lambda \\
&=
\int_I\left(g(\lambda),\mathcal F_0(\lambda)f\right)_{L^2(M_{\lambda})}
d\lambda, \quad
f \in \mathcal B.
\end{split}
\nonumber
\end{equation}
Therefore
$$
\left|\left(\int_I\mathcal F_0(\lambda)^{\ast}g(\lambda)d\lambda,f\right)
\right| \leq \|g\|\cdot\|\mathcal F_0f\| = \|g\|\cdot\|f\|.
$$
By Riesz' theorem, we then have
$$
\int_I\mathcal F_0(\lambda)^{\ast}g(\lambda)d\lambda \in L^2({\bf T}^d),
\quad
\|\int_I\mathcal F_0(\lambda)^{\ast}g(\lambda)d\lambda \| \leq \|g\|.
$$
Therefore for any compact interval $J \subset (0,d)\setminus{\bf Z}$,
$$
\|\int_J\mathcal F_0(\lambda)^{\ast}(\mathcal F_0f)(\lambda)d\lambda\| \leq
\|\mathcal F_0\chi_J(H_0)f\| = \|\chi_J(H_0)f\|,
$$
where $\chi_J$ is the characteristic function of $J$.
One can then show the existence of the strong limit
$$
\lim_{N\to\infty}\int_{I_N}\mathcal F_0(\lambda)^{\ast}(\mathcal F_0f)
(\lambda)d\lambda =:\tilde f.
$$
By Parseval's formula, we have for any $h \in L^2({\bf T}^d)$
\begin{eqnarray*}
(\tilde f,h) &=& \lim_{N\to\infty}\int_{I_N}(\mathcal F_0f(\lambda),
\mathcal F_0(\lambda)h)d\lambda \\
&=& (\mathcal F_0f,\mathcal F_0h) = (f,h),
\nonumber
\end{eqnarray*}
which implies $\tilde f = f$. \qed

\bigskip
Next let us construct the spectral representation for $H$.
We put
\begin{equation}
\mathcal F^{(\pm)}(\lambda) = \mathcal F_0(\lambda)\left(1 - VR(\lambda
\pm i0)\right), \quad
\lambda \in (0,d)\setminus\left({\bf Z}\cup\sigma_p(H)\right).
\nonumber
\end{equation}
Then by (\ref{S2LAP})
\begin{equation}
\mathcal F^{(\pm)}(\lambda) \in {\bf B}(\mathcal B\, ;L^2(M_{\lambda})).
\nonumber
\end{equation}


\begin{lemma}
For $\lambda \in (0,d)\setminus\left({\bf Z}\cup\sigma_p(H)\right)$, and
$f, g \in \mathcal B$
\begin{equation}
\left(\big(R(\lambda + i0) - R(\lambda - i0)\big)f,g\right) =
\left(\mathcal F^{(\pm)}(\lambda)f,\mathcal F^{(\pm)}
(\lambda)g\right)_{L^2(M_{\lambda})}.
\nonumber
\end{equation}
\end{lemma}

Proof.  We put
\begin{equation}
H_1 = H_0, \quad H_2 = H, \quad R_j(z) = (H_j - z)^{-1},
\nonumber
\end{equation}
\begin{equation}
G_{jk}(z) = (H_j - z)R_k(z),
\nonumber
\end{equation}
\begin{equation}
E_j'(\lambda) = \frac{1}{2\pi i}\left(R_j(\lambda + i0) -
R_j(\lambda - i0)\right).
\nonumber
\end{equation}
Then we have, by the resolvent equation,
\begin{equation}
\begin{split}
& \frac{1}{2\pi i}\left(R_k(\lambda + i\epsilon) -
R_k(\lambda - i\epsilon)\right) \\
=&\,
G_{jk}(\lambda \pm i\epsilon)^{\ast}\frac{1}{2\pi i}
\left(R_j(\lambda + i\epsilon) - R_j(\lambda - i\epsilon)
\right)G_{jk}(\lambda\pm i\epsilon).
\end{split}
\nonumber
\end{equation}
Letting $\epsilon \to 0$, we have for $f, g \in \mathcal B$
\begin{equation}
\left(E_k'(\lambda)f,g\right) = \left(E_j'(\lambda)G_{jk}
(\lambda \pm i0)f,G_{jk}(\lambda \pm i0)g\right).
\label{S2EkEjGjk}
\end{equation}
Let $j = 1$, $k = 2$. Since $\mathcal F^{(\pm)}(\lambda) =
\mathcal F_0(\lambda)G_{12}(\lambda \pm i0)$, the lemma follows
 if we replace $f,g$ in Lemma 3.1 by $G_{jk}(\lambda\pm i0)f$,
 $G_{jk}(\lambda\pm i0)g$. \qed

\bigskip
We define the operator $\mathcal F^{(\pm)}$  by $\big(\mathcal
F^{(\pm)}f\big)(\lambda) = \mathcal F^{(\pm)}(\lambda)f$ for $f
\in \mathcal B$. Let $\mathcal H_{ac}(H)$ be the absolutely
continuous subspace for $H$.


\begin{theorem}
(1) $\mathcal F^{(\pm)}$ is uniquely extended to a partial isometry with
 initial set $\mathcal H_{ac}(H)$ and final set $L^2((0,d);L^2(M_{\lambda});
 d\lambda)$. Moreover it diagonalizes $H$ :
\begin{equation}
\big(\mathcal F^{(\pm)}Hf\big)(\lambda) = \lambda \big(\mathcal
F^{(\pm)}f\big)(\lambda), \quad \forall f \in L^2({\bf T}^d).
\end{equation}
(2) The  following inversion formula holds:
\begin{equation}
f = \mathop{\rm s-lim}_{N\to\infty}\int_{I_N}
\mathcal F^{(\pm)}(\lambda)^{\ast}\big(\mathcal F^{(\pm)}f\big)
(\lambda)d\lambda, \quad \forall f \in {\mathcal H}_{ac}(H),
\end{equation}
where $I_N$ is a union of compact intervals $\subset (0,d)\setminus
\left({\bf Z}\cup\sigma_p(H)\right)$ such that $I_N \to (0,d)\setminus
\left({\bf Z}\cup\sigma_p(H)\right)$ as $N \to \infty$. \\
\noindent
(3) $\mathcal F^{(\pm)}(\lambda)^{\ast} \in {\bf B}(L^2(M_{\lambda});
\mathcal B^{\ast})$ is an eigenoperator for $H$ in the sense that
$$
(H - \lambda)\mathcal F^{(\pm)}(\lambda)^{\ast}\phi = 0, \quad
\forall \phi \in L^2(M_{\lambda}).
$$
(4) The wave operators
\begin{equation}
W^{(\pm)} = {\mathop{\rm s-lim}_{t\to\pm\infty}}\, e^{itH}e^{-itH_0}
\label{S2WaveOp}
\end{equation}
exist and are complete, i.e. the range of $W^{(\pm)}$ is equal to
$\mathcal H_{ac}(H)$. Moreover,
\begin{equation}
W^{(\pm)} = \big(\mathcal F^{(\pm)}\big)^{\ast}\mathcal F_0.
\label{S2StatTimedeWaOp}
\end{equation}
\end{theorem}

Proof. The proof of (1), (2), (3) is the same as that for Theorem 3.2
except for the surjectivity of $\mathcal F^{(\pm)}$. Since $V$ is trace class,
the existence and completeness of wave operators (\ref{S2WaveOp}) can be
proven by Rosenblum-Kato theory (see \cite{Ka76}, p. 542). The relation
(\ref{S2StatTimedeWaOp}) is also well-known, and we omit the proof
(see e.g. \cite{IsKu10}). We then have $\mathcal F^{(\pm)} =
\mathcal F_0\left(W^{(\pm)}\right)^{\ast}$. The completeness
of $W^{(\pm)}$ implies the surjectivity of $\mathcal F^{(\pm)}$. \qed


\subsection{Spectral representation on the lattice}

 Theorem 3.4 is transferred on the lattice space by $\mathcal U$.
We put $\widehat{\mathcal F}_0(\lambda)$, $\widehat{\mathcal F}^{(\pm)}
(\lambda)$, $\widehat{\mathcal F}_0$ and $\widehat{\mathcal F}^{(\pm)}$ by
\begin{equation}
\widehat{\mathcal F}_0(\lambda) = \mathcal F_0(\lambda) \, \mathcal U, \quad
\widehat{\mathcal F}^{(\pm)}(\lambda) = \mathcal F^{(\pm)}(\lambda) \,
\mathcal U,
\nonumber
\end{equation}
\begin{equation}
\widehat{\mathcal F}_0 = \mathcal F_0 \, \mathcal U, \quad
\widehat{\mathcal F}^{(\pm)} = \mathcal F^{(\pm)} \, \mathcal U.
\nonumber
\end{equation}
We also define
\begin{equation}
\widehat{\mathcal B} = \mathcal U^{-1}\mathcal B, \quad
\widehat{\mathcal B}^{\ast} = \mathcal U^{-1}{\mathcal B}^{\ast}.
\nonumber
\end{equation}


\begin{theorem}
(1) $\widehat{\mathcal F}^{(\pm)}$ is uniquely extended to a partial
isometry with intial set $\mathcal H_{ac}(\widehat H)$ and final set
$L^2({\bf T}^d)$. Moreover it diagonalizes $\widehat H$ :
\begin{equation}
\big(\widehat{\mathcal F}^{(\pm)}\widehat H\widehat f\big)(\lambda) =
\lambda \big(\widehat{\mathcal F}^{(\pm)}\widehat f\big)(\lambda).
\end{equation}
(2) The  following inversion formula holds:
\begin{equation}
\widehat f = \mathop{\rm s-lim}_{N\to\infty}\int_{I_N}
\widehat{\mathcal F}^{(\pm)}(\lambda)^{\ast}\left(\widehat
{\mathcal F}^{(\pm)}\widehat f\right)(\lambda)d\lambda, \quad
\forall \widehat f \in \mathcal H_{ac}(\widehat H),
\end{equation}
where $I_N$ is a union of compact intervals $\subset (0,d)
\setminus\left({\bf Z}\cup\sigma_p(\widehat H)\right)$ such
that $I_N \to (0,d)\setminus\left({\bf Z}\cup\sigma_p(\widehat H)
\right)$ as $N\to \infty$. \\
\noindent
(3) $\widehat{\mathcal F}^{(\pm)}(\lambda) \in {\bf B}(L^2(M_{\lambda})
\,;\,\widehat{\mathcal B}^{\ast})$ is an eigenoperator for $\widehat H$
in the sense that
$$
(\widehat H - \lambda)\widehat{\mathcal F}^{(\pm)}(\lambda)^{\ast}\phi = 0,
\quad \forall \phi \in L^2(M_{\lambda}).
$$
(4) The wave operators
\begin{equation}
\widehat W^{(\pm)} = {\mathop{\rm s-lim}_{t\to\pm\infty}}\,
e^{it\widehat H}e^{-it\widehat H_0}
\nonumber
\end{equation}
exist and are complete. Moreover,
\begin{equation}
{\widehat W}^{(\pm)} = \big(\widehat{\mathcal F}^{(\pm)}\big)^{\ast}
\widehat{\mathcal F}_0.
\nonumber
\end{equation}
\end{theorem}


\subsection{Generalized eigenvector}
It is customary to define the distribution $\delta(h(x) - \lambda) \in
\mathcal D'({\bf T}^d)$ by
\begin{equation}
\int_{{\bf T}^d}f(x)\delta(h(x)-\lambda)dx = \int_{M_{\lambda}}f(x)
\frac{dM_{\lambda}}{|\nabla_x h(x)|}, \quad f \in C({\bf T}^d).
\nonumber
\end{equation}
We then see that $\mathcal F_0(\lambda)^{\ast}$ defines a distribution
on ${\bf T}^d$ by the following formula
\begin{equation}
\mathcal F_0(\lambda)^{\ast}\phi = \phi(x)\delta(h(x) - \lambda).
\nonumber
\end{equation}
The right-hand side makes sense when, for example, $\phi \in C^{\infty}
(M_{\lambda})$ and is extended to a $C^{\infty}$-function near $M_{\lambda}$,
which is denoted by $\phi(x)$ again. The Fourier coefficients of
$\mathcal F_0(\lambda)^{\ast}\phi$ are then computed as
\begin{equation}
(2\pi)^{-d/2}\int_{{\bf T}^d}e^{in\cdot x}\delta(h(x)-\lambda)\phi(x)dx =
(2\pi)^{-d/2}\int_{M_{\lambda}}e^{in\cdot x}\phi(x)\frac{dM_{\lambda}}
{|\nabla_xh(x)|}.
\label{S3FourierCoeffi}
\end{equation}

We look for a parametrization of $M_{\lambda}$ suitable for the computation
in the next section.
Let us note
\begin{equation}
\frac{1}{2}\Big(d - \sum_{j=1}^d\cos x_j\Big) = \sum_{j=1}^d\sin^2
\left(\frac{x_j}{2}\right),
\nonumber
\end{equation}
which suggests that the variables $y = (y_1,\cdots,y_d) \in [-1,1]^d$:
\begin{equation}
y_j = \sin\frac{x_j}{2}, \quad x_j = 2\arcsin y_j
\nonumber
\end{equation}
are convenient to describe $H_0$. In fact, the map $x \to y$ is a
diffeomorphism between two tori:
$$
{\bf R}^d/(2\pi{\bf Z})^d = [-\pi,\pi]^d \ni x \to y \in [-1,1]^d=
{\bf R}^d/(2{\bf Z})^d.
$$
Conseqently,
\begin{equation}
x(\sqrt{\lambda}\theta) = \left(2\arcsin(\sqrt\lambda\theta_1),\cdots,
2\arcsin(\sqrt\lambda\theta_d)\right), \quad \theta \in S^{d-1},
\label{S2xsqrtlambadomega}
\end{equation}
gives a parameter representation of $M_{\lambda}$.
Passing to the polar coordinates $y = \sqrt{\lambda}\theta$, we also have
\begin{equation}
dx = J(y)dy = \frac{(\sqrt{\lambda})^{d-2}J(\sqrt{\lambda}\theta)}{2}
d\lambda d\theta,
\label{S3Jydy}
\end{equation}
which implies
\begin{equation}
\frac{dM_{\lambda}}{|\nabla_xh(x)|} = \frac{(\sqrt{\lambda})^{d-2}
J(\sqrt{\lambda}\theta)}{2}d\theta.
\label{S3dMnablap}
\end{equation}

We define $\widehat \psi^{(0)}(n,\lambda,\theta)$ by
\begin{equation}
\begin{split}
\widehat\psi^{(0)}(n,\lambda,\theta)
&= (2\pi)^{-d/2}\frac{(\sqrt\lambda)^{d-2}}{2}\,
e^{in\cdot x(\sqrt\lambda\theta)}\, J(\sqrt\lambda\theta) \\
&= (2\pi)^{-d/2}2^{d-1}(\sqrt\lambda)^{d-2}\, e^{in\cdot
x(\sqrt{\lambda}\theta)}\,
\frac{\chi(\sqrt\lambda\theta)}{{\mathop\Pi_{j=1}^d}\cos
\big(x_j(\sqrt{\lambda}\theta)/2\big)}.
\end{split}
\label{S3psi0nlambdatheta}
\end{equation}
where $\chi(y)$ is the characteristic function of $[-1,1]^d$.
By (\ref{S3FourierCoeffi}) and (\ref{S3dMnablap}), we have for
$\phi \in L^2(M_{\lambda})$
\begin{equation}
\big(\widehat{\mathcal F}_0(\lambda)^{\ast}\phi\big)(n) =
\int_{S^{d-1}}{\widehat\psi}^{(0)}(n,\lambda,\theta)
\phi(x(\sqrt{\lambda}\theta))d\theta.
\label{S3Formulahatfast}
\end{equation}
One can also see that if $\widehat f$ is compactly supported
\begin{equation}
\big(\widehat {\mathcal F}_0(\lambda)\widehat f\big)
(x(\sqrt{\lambda}\theta))= (2\pi)^{-d/2}\sum_{n\in{\bf Z}^d}
e^{-in\cdot x(\sqrt{\lambda}\theta)}\widehat f(n).
\label{S3Formulahatf}
\end{equation}


\subsection{Scattering matrix}
 The scattering operator $\widehat S$ is defined by
\begin{equation}
\widehat S = \big(\widehat W_{+}\big)^{\ast}\widehat W_{-}.
\nonumber
\end{equation}
We conjugate it by the spectral representation. Let
\begin{equation}
\mathcal S = \widehat {\mathcal F}_0\, \widehat S\, \big(\widehat{\mathcal F}
_0\big)^{\ast},
\nonumber
\end{equation}
which is unitary on $L^2((0,d);L^2(M_{\lambda});d\lambda)$.
Since $\mathcal S$ commutes with $\widehat H_{0}$,
${\mathcal S}$ is written as a direct integral
\begin{equation}
{\mathcal S} = \int_{(0,d)}{\oplus}{\mathcal S}(\lambda)d\lambda.
\nonumber
\end{equation}
The S-matrix, ${\mathcal S}(\lambda)$, is unitary on $L^2(M_{\lambda})$
and has the
following representation.


\begin{theorem}
Let $\lambda \in (0,d)\setminus\big({\bf Z}\cup\sigma_p(H)\big)$.
Then ${\mathcal S}(\lambda)$ is written as
\begin{equation}
{\mathcal S}(\lambda) = 1 - 2\pi i A(\lambda),
\nonumber
\end{equation}
\begin{equation}
A(\lambda) = \widehat{\mathcal F}_0(\lambda)\widehat V
\widehat{\mathcal F}_0(\lambda)^{\ast}-
\widehat{\mathcal F}_0(\lambda)\widehat V\widehat
R(\lambda + i0)\widehat V\widehat{\mathcal F}_0(\lambda)^{\ast}.
\label{S3Alambda}
\end{equation}
\end{theorem}

Since the proof is well-known, we omit it (see e.g. \cite{IsKu10}).

\bigskip

By (\ref{S3psi0nlambdatheta}) and (\ref{S3Formulahatf}), the first term
of the right-hand side of (\ref{S3Alambda}) has an integral kernel
\begin{equation}
(2\pi)^{-d}\frac{\lambda^{d-2}}{4}J(\sqrt\lambda\theta)J(\sqrt\lambda\theta')
\sum_{n \in {\bf Z}^d}e^{-in\cdot(x(\sqrt{\lambda}\theta) -
x(\sqrt{\lambda}\theta'))}\widehat V(n).
\label{S3SmatrixBorn}
\end{equation}
The  2nd term of the right-hand side of (\ref{S3Alambda}) has the
following kernel
\begin{equation}
- (2\pi)^{-d/2}\frac{(\sqrt{\lambda})^{d-2}}{2}J(\sqrt\lambda\theta)
\sum_{n\in{\bf Z}^d}e^{-in\cdot x(\sqrt{\lambda}\theta)}\widehat V(n)
\left(\widehat R(\lambda + i0)\widehat V\widehat\psi^{(0)}(\lambda,\theta')
\right)(n).
\label{S3SmatrixRemainder}
\end{equation}


\section{Inverse scattering}

In this section we  prove that the potential $\widehat V$ is uniquely
reconstructed from the scattering matrix for all energies.
We first consider the analytic continuation of $x(\sqrt{\lambda}\theta)$
defined by (\ref{S2xsqrtlambadomega}).


\begin{lemma}
Let $\tau$ be a constant such that  $-1 < \tau < 1$, $\tau \neq 0$.
Then $f(z,\tau) = 2\arcsin(z\tau)$ is analytic with respect to
$z \in {\bf C}\setminus\Big((-\infty,-1/|\tau|]\cup[1/|\tau|,\infty)\Big)$.
 Moreover, letting  $\epsilon(\tau) = \tau/|\tau|$, we have as $N \to \infty$,
\begin{equation}
\begin{split}
{\rm Re}\, f(N+ i,\tau) & \equiv \epsilon(\tau)\Big(\pi - \frac{2}{N} +
O(\frac{1}{N^3})\Big), \quad {\rm mod} \quad 2\pi{\bf Z} \\
{\rm Im}\,f(N+i,\tau) & = \epsilon(\tau)\Big(2\log N + \log(4\tau^2) +
O(\frac{1}{N^2})\Big).
\end{split}
\nonumber
\end{equation}

\end{lemma}

Proof. We take the branch of $\arcsin(z) = u(z) + iv(z)$ so that it
is single-valued analytic on ${\bf C}\setminus\big((-\infty,-1]
\cup[1,\infty)\big)$, and $0 < u(s) < \pi/2$ when $0 < s < 1$. Then we have
\begin{equation}
 \sin(u)\cosh(v) = {\rm Re}\, z, \quad
 \cos(u)\sinh(v) =  {\rm Im}\, z.
\label{S4FRealImagi}
\end{equation}

Let us note that
\begin{equation}
\cos(u) > 0 \quad {\rm and} \quad \pm \sinh(v) > 0, \quad
{\rm if} \quad \pm {\rm Im}\, z > 0.
\label{S4cos>0}
\end{equation}
In fact, by the 2nd equation of (\ref{S4FRealImagi}), $\cos(u(z))$ and
$\sinh(v(z))$ do not vanish if ${\rm Im}\,z \neq 0$, and $\cos(u(s)) > 0$,
 $v(s) = 0$ when $s \in (0,1)$. So, $\cos(u(z)) > 0$ when ${\rm Im}\, z
 \neq 0$, and again by the 2nd equation of (\ref{S4FRealImagi}),
 $\sinh(v(z))$ and $ {\rm Im}\, z > 0$
have the same sign.

Let $u_N = u((N+i)\tau), v_N = v((N+i)\tau)$. Then by (\ref{S4FRealImagi})
we have $\sinh(v_N) =
\tau/\cos(u_N)$.
Plugging this with
$$
\left(1 - \cos^2(u_N)\right)\left(1 + \sinh^2(v_N)\right) = N^2\tau^2,
$$
and letting $t_N = \cos^2(u_N)$, we get the equation
$$
t_N^2 + (N^2\tau^2 + \tau^2 - 1)t_N - \tau^2 = 0.
$$
Since $t_N > 0$, by solving this equation, we have
$$
t_N = N^{-2} + O(N^{-4}).
$$
Since $\cos(u_N) > 0$, we have
\begin{equation}
\cos(u_N) = N^{-1} + O(N^{-3}).
\label{S4cosuN}
\end{equation}
This, combined with (\ref{S4FRealImagi}) for $z = (N+i)\tau$, then yields
\begin{equation}
\sinh(v_N) = \tau N + O(N^{-1}).
\label{S4sinhvN}
\end{equation}

If $\tau > 0$, then $v_N > 0$, and by (\ref{S4FRealImagi}) with
$z = (N+i)\tau$, $\sin(u_N) > 0$.
From (\ref{S4cosuN}), we then get
\begin{equation}
u_N \equiv \frac{\pi}{2}- N^{-1} + O(N^{-3}) \quad \rm{mod} \quad
2\pi{\bf Z}.
\label{S4uN}
\end{equation}
From (\ref{S4sinhvN}), we have
$$
e^{2v_N} - 2\left(\tau N + O(N^{-1})\right)e^{v_N} - 1 = 0,
$$
which implies
\begin{equation}
v_N = \log(2\tau N) + O(N^{-2}).
\label{S4vN}
\end{equation}
Since $f(z,\tau) = 2\arcsin(z\tau)$ and $\arcsin(-z) = - \arcsin(z)$,
(\ref{S4uN}jand (\ref{S4vN}) prove the lemma.
\qed

\medskip
We define the $l^1$-norm of $n = (n_1,\cdots,n_d) \in {\bf Z}^d$ by
\begin{equation}
|n|_{l^1} = \sum_{j=1}^d|n_j|.
\end{equation}

We also introduce the following  notation. For $n \in {\bf Z}^d$, we
define $\widehat B(n) \in {\bf B}(l^2({\bf Z}^d);{\bf C})$ and $\widehat K(n)
 \in {\bf B}({\bf C};l^2({\bf Z}^d))$ by
\begin{equation}
\widehat B(n)\widehat f = \widehat f(n),
\label{S4hatIm}
\end{equation}
\begin{equation}
\big(\widehat K(n)c\big)(m) = c\,\delta_{mn}.
\label{S4widehatJm}
\end{equation}
Note that
\begin{equation}
\widehat P(n) = \widehat K(n)\widehat B(n).
\label{S4P=KB}
\end{equation}

Let $r_0(k,z)$ be defined by
\begin{equation}
r_0(k,z) = (2\pi)^{-d}\int_{T^d}\frac{e^{ik\cdot x}}{h(x) - z}dx.
\label{S3r0kz}
\end{equation}
Then the resolvent $\widehat R_0(z) = (\widehat H_0 -z)^{-1}$ is written as
\begin{equation}
\big(\widehat R_0(z)\widehat f\big)(m) =
\sum_{n\in{\bf Z}^d}r_0(m-n,z)\widehat f(n).
\label{S3hatR0(z)}
\end{equation}


\begin{lemma}
For any $m, n$ and $|z| > d$
\begin{equation}
\widehat B(m)\widehat R_0(z) \widehat
K(n) = \sum_{s=|m-n|}^{\infty}c_s(m-n)\, z^{-s-1},
\label{S4R0(z)Expand}
\end{equation}
where $c_s(k)$ is a constant satisfying
\begin{equation}
|c_s(k)| \leq d^s, \quad s = 0,1,2,\cdots.
\label{S4Rosbound}
\end{equation}
In particular,  for $|z| > 2d$
\begin{equation}
\|\widehat B(m)\widehat R_0(z) \widehat K(n)\| \leq |z|^{-1-|m-n|_{l^1}}.
\label{S4R0zdecay}
\end{equation}
\end{lemma}

Proof. Using $(h(x) - z)^{-1} = - \sum_{s=0}^{\infty}z^{-s-1}h(x)^s$,
we have for large $|z|$
\begin{equation}
r_0(k,z) = \sum_{s=0}^{\infty}z^{-s-1}c_s(k),
\label{S4r0kzexpand}
\end{equation}
$$
c_s(k) = - (2\pi)^{-d}\int_{T^d}h(x)^se^{ik\cdot x}dx.
$$
Since $|h(x)| \leq d$, we have
\begin{equation}
|c_s(k)| \leq d^s.
\label{S4Cskboubd}
\end{equation}
 Hence the series (\ref{S4r0kzexpand}) is absolutely convergent for $|z| > d$.
Note that $h(x)^s$ is a sum of terms of the form
$$
\left(\cos x_1\right)^{\alpha_1}\cdots\left(\cos x_d\right)^{\alpha_d}, \quad
0 \leq \alpha_1 + \cdots + \alpha_d \leq s, \quad 0 \leq \alpha_j \leq s.
$$
If $s < |k|_{l^1} $, we have $\alpha_j < |k_j| $ for some $j$, which implies
$$
\int_{T^d}h(x)^se^{ik\cdot x}dx = 0, \quad {\rm if} \quad s < |k|_{l^1}.
$$
Then we have
\begin{equation}
c_s(k) = 0, \quad {\rm if} \quad s < |k|_{l^1}.
\label{S4Csk}
\end{equation}
By (\ref{S3hatR0(z)}) and (\ref{S4r0kzexpand}), we have
$$
\widehat B(m)\widehat R_0(z)\widehat K(n) = \sum_{s=0}^{\infty}z^{-s-1}c_s(m-n).
$$
Then lemma then follows from (\ref{S4Cskboubd}) and (\ref{S4Csk}). \qed


\begin{lemma}
For any $m, n$, there exists a constant $C_{mn}$ such that if $|z| >
\|\widehat H\|+1$,
$$
\|\widehat B(m)\widehat R(z)\widehat K(n)\| \leq
C_{mn}(1 + |z|)^{-1-|m-n|_{l^1}}.
$$
\end{lemma}

Proof. Let $p = |m-n|_{l^1}$. By the perturbation expansion, we have
$$
\widehat R(z) = \widehat R_0(z) - \widehat R_0(z)\widehat V\widehat R_0(z)
+ \cdots + (-1)^p\widehat R_0(z)\widehat V\cdots\widehat V
\widehat R_0(z) + O(z^{-p-2}).
$$
Multiply this equality by $\widehat B(m)$ and $\widehat K(n)$.
Then by Lemma 4.2, the first term decays like $O(|z|^{-p-1})$.
Next we look at the term
$$
\widehat B(m)\widehat R_{0}(z)\widehat V\cdots\widehat V
\widehat R_{0}(z)\widehat K(n),
$$
consisting of $k$ numbers of $\widehat V$ and $k+1$ numbers of
$\widehat R_0(z)$, where $1 \leq k \leq p$. It is rewritten as a
finite linear combination of terms
\begin{equation}
\widehat B(m)\widehat R_0(z)\widehat P(r^{(1)})\widehat R_0(z)
\widehat P(r^{(2)})\cdots
\widehat P(r^{(k)})\widehat R_0(z)\widehat K(n),
\label{S4pmr0...}
\end{equation}
since $\widehat V = \sum_{|r| < c}\widehat V(r)\widehat P(r)$ for some $c > 0$.
We put
$$
\epsilon_1 = |m-r^{(1)}|_{l^1}, \quad \epsilon_2 = |r^{(1)}-r^{(2)}|_{l^1},
\quad \cdots, \epsilon_{k+1} = |r^{(k)}-n|_{l^1}.
$$
By (\ref{S4P=KB}), (\ref{S4pmr0...}) is written as the product
$$
\widehat B(m)\widehat R_0(z)\widehat K(r^{(1)})\cdot
\widehat B(r^{(1)})\widehat R_0(z)\widehat K(r^{(2)})\cdot\cdots
\cdot\widehat B(r^{(k)})\widehat R_0(z)\widehat K(r^{(n)}).
$$
By (\ref{S4R0zdecay}), this decays like $|z|^{-(k+1 + \epsilon_1 + \cdots + \epsilon_{k+1})}$.
Since $|m-n|_{l^1} = p$, we have
$$
\epsilon_1 + \cdots + \epsilon_{k+1} \geq p.
$$
Taking notice of $k+1 + \epsilon_1 + \cdots + \epsilon_{k+1} \geq 2 + p$, we have proven the lemma. \qed

\bigskip
We can now solve the inverse problem for $\widehat H$.

\begin{theorem}
Suppose $\widehat V$ is compactly supported. Then
from the scattering amplitude $A(\lambda)$ for all $\lambda \in (0,d)\setminus\left({\bf Z}\cup\sigma_p(H)\right)$, one can reconstruct $\widehat V$.
\end{theorem}

Proof.
Let $A(\lambda;\theta,\theta)$ be the integral kernel of $A(\lambda)$. Let $\sqrt\lambda = k$, and put
\begin{equation}
B(k;\theta,\theta') = \frac{4(2\pi)^d}{k^{2d-4}}\left(J(k\theta)J(k\theta')\right)^{-1}A(k^2;\theta,\theta').
\nonumber
\end{equation}
In view of (\ref{S3psi0nlambdatheta}), (\ref{S3SmatrixBorn}) and (\ref{S3SmatrixRemainder}), we can rewrite it as
\begin{equation}
B(k;\theta,\theta') = B_0(k;\theta,\theta') - B_1(k;\theta,\theta'),
\label{S3Bktheta}
\end{equation}
\begin{equation}
B_0(k;\theta,\theta') = \sum_{n\in{\bf Z}^d}e^{-in\cdot(x(k\theta) - x(k\theta'))}\widehat V(n),
\label{S3B0ktheta}
\end{equation}
\begin{equation}
B_1(k;\theta,\theta') = - \sum_{n\in{\bf Z}^d}e^{-in\cdot x(k\theta)}\widehat V(n)
\left(\widehat R(k^2+i0)\widehat V\widehat\varphi^{(0)}(k,\theta')\right)(n),
\label{S3B1ktheta}
\end{equation}
\begin{equation}
\widehat\varphi^{(0)}(k,\theta') = \left(e^{in\cdot x(k\theta')}\right)_{n\in{\bf Z}^d}.
\label{S3hatvarphi0}
\end{equation}
Let $\theta_j \neq 0, \ \theta_j' \neq 0$ $(j = 1,\cdots, d)$, and put $\zeta(z,\theta) = (f(z,\theta_1),\cdots,f(z,\theta_d))$, where $f(z,\tau)$ is defined in Lemma 4.1. Then $B_0(k;\theta,\theta')$ and $B_1(k;\theta,\theta')$ have analytic continuations $B_0(z;\theta,\theta')$ and $B_1(z;\theta,\theta')$, which are defined with $k$ replaced by $z$ in the upper-half plane and $x(k\theta)$ by $\zeta(z,\theta)$.
We put
$$
S(n) = n_1 + \cdots + n_d, \quad n \in {\bf Z^d}.
$$

We now take $\theta_j > 0$ and $\theta_j' < 0$, $(j = 1,\cdots, d)$. Then by Lemma 4.1, as $z = N + i \to \infty$,
\begin{equation}
e^{-in\cdot\zeta(z,\theta)} = e^{-iS(n)\pi}(4N^2)^{S(n)}{\mathop\Pi_{j=1}^d}
(\theta_j)^{2n_j}(1 + O(N^{-1})),
\label{S4einzetaexpand}
\end{equation}
\begin{equation}
e^{in\cdot\zeta(z,\theta')} = e^{iS(n)\pi}(4N^2)^{S(n)}{\mathop\Pi_{j=1}^d}
\big(\theta_j'\big)^{2n_j}(1 + O(N^{-1})).
\label{S4e-inzetaexpand}
\end{equation}

The recontsruction procedure for $\widehat V(n)$ goes inductively with respect to $S(n)$.  Let us assume that
\begin{equation}
{\rm supp}\,\widehat V \subset \{(n_1,\cdots,n_d) \, ; \, |n_j| \leq M\}.
\label{S4SupportVbound}
\end{equation}
We use (\ref{S3B0ktheta}), (\ref{S4einzetaexpand}) and (\ref{S4e-inzetaexpand}) to compute the asymptotic expansion of $B_0(z_N;\theta,\theta')$, $z_N = N + i$,  as $N \to \infty$. Then the largest contribution arises from the term  for which $S(n)$ is the largest, i.e. $n = (M,\cdots,M)$. Therefore, we have
\begin{equation}
B_0(z_N\, ;\, \theta,\theta') \sim (2N)^{4dM}{\mathop\Pi_{j=1}^d}\big(\theta_j\theta_j'\big)^{2M}\widehat V(n^{(M)}),
\label{S4B0topterm}
\end{equation}
 where $n^{(M)} = (M,\cdots,M)$. Since
\begin{equation}
\|\hat R((N+i)^2)\| = O(N^{-2}),
\label{S4RN+i2decay}
\end{equation}
 using (\ref{S3B1ktheta}), (\ref{S4einzetaexpand}) and (\ref{S4e-inzetaexpand}), we see that
\begin{equation}
B_1(z_N\, ;\, \theta,\theta') = O(N^{4dM-2}).
\label{S4B1topterm}
\end{equation}
By (\ref{S4B0topterm}) and (\ref{S4B1topterm}), we can compute $\widehat V(n^{(M)})$ from the asymptotic expansion of $B(z_N\, ;\, \theta,\theta')$.

Assume that we have computed $\widehat V(n)$ for $S(n) > p$. Then
\begin{equation}
\begin{split}
B_0(z_N\, ;\, \theta,\theta') &- \sum_{S(n)>p}e^{-in\cdot(\zeta(z,\theta)-\zeta(z,\theta'))}\widehat V(n) \\
& \sim (2N)^{4p}\sum_{S(n)=p}{\mathop\Pi_{j=1}^d}(\theta_j\theta_j')^{2n_j}\widehat V(n).
\end{split}
\label{S4B0expand}
\end{equation}
The image of the map
$$
S_-^{d-1}\times S_+^{d-1} \ni (\theta,\theta') \to (\theta_1\theta_1',\cdots,\theta_d\theta_d')
$$
contains an open set in ${\bf R}^d$, where
$$
S_+^{d-1} = \{\theta \in S^{d-1}\,;\,\theta_j > 0, \ \forall j\}, \quad
S_-^{d-1} = \{\theta \in S^{d-1}\,;\,\theta_j < 0, \ \forall j\}.
$$
Therefore, one can compute $\widehat V(n)$ for $S(n) = p$ from (\ref{S4B0expand}).

We show that $B_1(z_N;\theta,\theta') = O(N^{4p-2})$ up to terms which are already known. We rewrite $B_1(z_N;\theta,\theta')$ as
\begin{equation}
\begin{split}
B_1(z_N;\theta,\theta') &= \sum_{m,n}e^{-im\cdot\zeta(z_N,\theta)}e^{in\cdot\zeta(z_N,\theta')}\cdot\widehat V(m)\widehat V(n)\cdot
\widehat B(m)\widehat R(z_N^2)\widehat K(n).
\end{split}
\nonumber
\end{equation}
We split this sum into 4 parts:
$$
\sum_{S(m),S(n)>p} + \sum_{S(m)>p,S(n)\leq p} + \sum_{S(m)\leq p,S(n)>p} + \sum_{S(m),S(n)\leq p} =: I_1 + I_2 + I_3 + I_4.
$$
Note that by (\ref{S4einzetaexpand}) and (\ref{S4e-inzetaexpand}), we have
\begin{equation}
e^{-im\cdot\zeta(z_N,\theta)}e^{in\cdot\zeta(z_N,\theta')} = O(N^{2(S(m)+S(n))}).
\label{S4eimein}
\end{equation}
Then by (\ref{S4RN+i2decay}) and (\ref{S4eimein}), we have
\begin{equation}
I_4 = O(N^{4p-2}).
\label{S3I4decay}
\end{equation}
By Lemma 4.3 and (\ref{S4eimein}), $I_3 = O(N^{2(S(m)+ S(n) -1 -|m-n|_{l^1})})$.
Since
$$
S(n) - S(m) = \sum_{j=1}^d(n_j - m_j) \leq |m - n|_{l^1},
$$
we have
\begin{equation}
\begin{split}
& 2S(m) + 2S(n) - 2|m-n|_{l^1} - 2 \\
& = 4S(m) + 2(S(n) - S(m) - |m-n|_{l^1}) - 2 \leq 4p - 2,
\end{split}
\nonumber
\end{equation}
which proves
\begin{equation}
I_3 = O(N^{4p-2}).
\label{S4I3decay}
\end{equation}
Similarly, we can prove
\begin{equation}
I_2 = O(N^{4p-2}).
\label{S4I2decay}
\end{equation}

We finally observe $I_1$. We put
$$
\widehat V_{\leq p} = \sum_{S(n)\leq p}\widehat V(n)\widehat P(n), \quad
\widehat V_{> p} = \sum_{S(n)> p}\widehat V(n)\widehat P(n),
$$
$$
\widehat H_{> p} = \widehat H_0 + \widehat V_{>p}, \quad
\widehat R_{>p}(z) = (\widehat H_{>p} - z)^{-1}.
$$
By the resolvent equation, $I_1$ is split into 2 parts
\begin{equation}
\begin{split}
I_1 &= \sum_{S(m),S(n)>p}e^{-im\cdot\zeta(z_N,\theta)}e^{in\cdot\zeta(z_N,\theta')}\cdot \widehat V(m)\widehat V(n)\cdot\widehat B(m)
\widehat R_{>p}(z_N^2)\widehat K(n)  \\
&- \sum_{S(m),S(n)>p}e^{-im\cdot\zeta(z_N,\theta)}e^{in\cdot\zeta(z_N,\theta')}\cdot \widehat V(m)\widehat V(n)\cdot\widehat B(m)
\widehat R_{>p}(z_N^2)\widehat V_{\leq p}\widehat R(z_N^2)\widehat K(n).
\end{split}
\nonumber
\end{equation}
The 1st term of the right-hand side is a known term, since we have already reconstructed $\widehat V(n)$ for $S(n) > p$. The 2nd term is a linear combination of terms
\begin{equation}
e^{-im\cdot\zeta(z_N,\theta)}e^{in\cdot\zeta(z_N,\theta')}\cdot\widehat B(m)\widehat R_{>p}(z_N^2)\widehat P(k)\widehat R(z_N^2)\widehat K(n),
\label{S4eimRpetc}
\end{equation}
where $S(k) \leq p$.
By Lemma 4.3, it decays like $O(N^{2(S(m)+S(n)-|m-k|_{l^1}-|n-k|_{l^1}-2)})$.
Using
\begin{equation}
\begin{split}
S(m) - |m-k|_{l^1} &= \sum_j(m_j - k_j) - \sum_j|m_j-k_j| + \sum_j k_j \\
& \leq S(k) \leq p,
\end{split}
\nonumber
\end{equation}
we see that (\ref{S4eimRpetc}) decays like $O(N^{4p-4})$. Therefore, we have
\begin{equation}
I_1 = O(N^{4p-2}),
\label{S4I1decay}
\end{equation}
up to known terms. By virtute of (\ref{S3I4decay}) $\sim$ (\ref{S4I1decay}),
we have completed the proof of the theorem. \qed


\section {Estimates of the free resolvent near the critical values}

The purpose of this section is to derive  estimates of the resolvent $\widehat R_0(z) = (\widehat H_0-z)^{-1}$ in weighted Hilbert spaces. Equivalently, we consider the operator norm on $L^2({\bf Z}^d)$ of $\widehat q\,\widehat R_0(z)\,\widehat q$, where $\widehat q$ is the operator of multiplication by $\widehat q \in \ell^{\infty}({\bf Z}^d)$ : $(\widehat q\,\widehat f)(n) = \widehat q(n)\widehat f(n)$. In particular, $\widehat \rho_s$ is the
operator defined by
$$
\widehat \rho_s(n) = (1 + |n|^2)^{-s/2}, \quad s \in {\bf R}.
$$

We put
\begin{equation}
{\bf D}_1=\{z \in {\bf C}\, ;\, 0 < |z|< 1\},
\nonumber
\end{equation}
\begin{equation}
\Lambda_d = \rho(\widehat H_0) = {\bf C}\setminus  [0,d],
\nonumber
\end{equation}
\begin{equation}
\lambda(z) = \frac{1}{4}\left(2 -z -  \frac{1}{z}\right).
\label{S5lambda(z)}
\end{equation}
For $w = re^{i\theta}$ with $0 < \theta < 2\pi$, we take the branch $\sqrt{w} = \sqrt{r}e^{i\theta/2}$. For $a, b \in {\bf C}$, we put
\begin{equation}
I(a,b) = \{ta + (1- t)b \, ; \, 0 \leq t \leq 1\}.
\nonumber
\end{equation}


\begin{lemma}
$\lambda(z)$ is a conformal map from ${\bf D}_1$ onto $\Lambda_1$, and its inverse is given by
\begin{equation}
z(\lambda) = 2\lambda - 1 - 2\sqrt{\lambda(\lambda -1)}.
\label{S5z(lambda)}
\end{equation}
\end{lemma}

Proof. Since the map
$$
\{0 < |z| < 1\} \ni z  \to \frac{1}{2}\left(z+\frac{1}{z}\right) \in {\bf C}\setminus[-1,1]
$$
is conformal, so is $\lambda(z)$ from ${\bf D}_1$ onto $\Lambda_1$.
By solving the equation $z^2 + (4\lambda-2)z + 1=0$, we have the inverse map
$z = 2\lambda - 1 \pm 2\sqrt{\lambda(\lambda-1)}$.
For $\lambda > 1$, $|2\lambda - 1 - 2\sqrt{\lambda(\lambda-1)}| < 1$. Therefore, we obtain (\ref{S5z(lambda)}). \qed

\medskip
The following lemma is proved by Lemma 2.4  and Theorem 2.6 (2). Note that $\widehat H_0$ has no eigenvalues.


\begin{lemma}
 For $s > 1/2$, the operator-valued function
$\widehat \rho_s\, \widehat R_0(\lambda)\, \widehat \rho_s$ is analytic with respect to $\lambda \in \Lambda_d$, and has continuous boundary values when $\lambda$ approaches $E \pm i0$, $E \in (0,d)\setminus{\bf Z}$.
\end{lemma}

We study estimates for $\widehat R_0(E \pm i0)$ when $E \in \sigma(\widehat H_0) = [0,d]$ is close to $\{0,1,\cdots,d\}$, the set of critical values of $h(x)$.
Let us begin with the case $d = 1$.


\begin{lemma} Assume $d = 1$.\\
\noindent
(1)  $r_0(n,z)$ defined by (\ref{S3r0kz}) has the following representation
\begin{equation}
 r_0(n,\lambda(z)) = \frac{4z^{|n|}}{z - 1/z}
=-\frac{z^{|n|}}{\sqrt{\lambda(z)(\lambda(z)-1)}},\quad {\rm for} \quad
(n,z)\in {\bf Z}\times {\bf D}_1.
\nonumber
\end{equation}
Moreover, $r_0(n,\lambda(z))$ has a
meromorphic continuation from ${\bf D}_1$ into
${\bf C}$.\\
\noindent
(2) Let $\|\cdot\|_{HS}$ be the Hilbert-Schmidt norm on $\ell^2({\bf Z})$,
and take $\widehat q_j=(\widehat q_j(n))_{n\in {\bf Z}}\in \ell^2({\bf Z}),
\ j= 1,2$. Then we have for $\lambda \in \Lambda_1$
\begin{equation}
\|\widehat q_1\, \widehat R_0(\lambda)\, \widehat q_2\|_{HS}\le
\frac{\|\widehat q_1\|_{\ell^2({\bf Z})}\|\widehat q_2\|_{\ell^2({\bf Z})}}
{|\sqrt{\lambda(\lambda-1)}|}.
\nonumber
\end{equation}
(3) Let 
$u(\nu)=\sqrt{\nu(\nu-1)}$, and take $\lambda, \lambda \in \Lambda_1$ such that $I(\lambda,\lambda_1) \subset \Lambda_1$. We put
\begin{equation}
\begin{split}
M(\lambda,\lambda_1)& =\left(1 + 2\max_{\nu\in I(\lambda,\lambda_1)}|u'(\nu)|\right)
\left(1 + \frac{1}{|u(\lambda)|} + \frac{1}{|u(\lambda_1)|}\right), \\
N(\lambda,\lambda_1) & = \left|\frac{1}{u(\lambda)} - \frac{1}{u(\lambda_1)}\right|.
\end{split}
\nonumber
\end{equation}
Then for any $0 \leq \alpha \leq 1$, we have the following pointwise H{\"o}lder estimate
\begin{equation}
|r_0(n,\lambda)-r_0(n,\lambda_1)|\le |\lambda - \lambda_1|^{\alpha}
(1 + |n|)^{\alpha}M(\lambda,\lambda_1)^{\alpha}N(\lambda,\lambda_1)^{1-\alpha},
\label{S5r0lambda-r0lambda1}
\end{equation}
and the following H{\"o}lder estimate of the Hilbert-Schmidt norm
\begin{equation}
\label{holder1} \|\widehat\rho_{\alpha}\widehat q_1\, (\widehat
R_0(\lambda)-\widehat R_0(\lambda_1))\,\widehat\rho_{\alpha}\widehat q_2\|_{HS}\le
|\lambda-\lambda_1|^{\alpha}C_{\alpha}(\lambda,\lambda_1) \|\widehat q_1\|_{\ell^2(\bf
Z)}\|\widehat q_2\|_{\ell^2(\bf Z)},
\end{equation}
\begin{equation}
\label{Cholder1} C_{\alpha}(\lambda,\lambda_1)= M(\lambda,\lambda_1)^{\alpha}
N(\lambda,\lambda_1)^{1-\alpha}.
\end{equation}
In particular, there exists a constant $C_{\alpha}' > 0$ such that
\begin{equation}
\label{holder3} 
\begin{split}
& \|\widehat\rho_{\alpha}\widehat q_1\, (\widehat
R_0(\lambda)-\widehat R_0(\lambda_1))\,\widehat\rho_{\alpha}\widehat q_2\|_{HS} \\
& \le
\frac{C_{\alpha}'|\lambda-\lambda_1|^{\alpha}}{|\lambda(\lambda -1)\lambda_1(\lambda_1 - 1)|^{(1+3\alpha)/2}}\, \|\widehat q_1\|_{\ell^2(\bf
Z)}\|\widehat q_2\|_{\ell^2(\bf Z)},
\end{split}
\end{equation}
if $|\lambda|, |\lambda_1|\leq 2$, and $I(\lambda,\lambda_1) \subset \Lambda_1$.
\end{lemma}

Proof. To prove (1), we first note by residue calculus
$$
r_0(n,\lambda(z))=\frac{1}{\pi}\int_0^{2\pi}\frac{e^{-inx}dx}{1-\cos x-2\lambda(z)} =\frac{2}{\pi
i}\int_{|w|=1}\frac{w^{-n}dw}{(w-z)(w- 1/z)} =
\frac{4z^{|n|}}{z - 1/z}.
$$
By (\ref{S5z(lambda)}), we have $1/z = 2\lambda-1+2\sqrt{\lambda(\lambda-1)}$.
 Hence
$z - 1/z = - 4\sqrt{\lambda(\lambda-1)}$, which proves (1).

Using (\ref{S3hatR0(z)}) and (1), we obtain for $z \in {\bf D}_1$
$$
\|\widehat q_1\, \widehat R_0(\lambda(z))\,\widehat q_2\|_{HS}^2=
\sum_{n,m\in {\bf Z}}
\frac{|z|^{2|n-m|}}{|\lambda(z)(\lambda(z)-1)|}|\widehat q_1(n)|^2|
\widehat q_2(m)|^2 \\
$$
$$
\leq
\frac{\|\widehat q_1\|^2\|\widehat q_2\|^2}{|\lambda(z)(\lambda(z)-1)|},
$$
where $\|\widehat q_j\|^2=\sum_{n\in {\bf Z}} |\widehat q_j(n)|^2$.
This proves (2).

We prove (3). Let $\lambda_1,\lambda\in {\bf C_+}$ and $z = z(\lambda), z_1=z(\lambda_1)$. Then (1) and
(\ref{S5z(lambda)}) give
$$
r_0(n,\lambda)-r_0(n,\lambda_1) =-\frac{z^{|n|}}{u(\lambda)}+
\frac{z_1^{|n|}}{u(\lambda_1)}=\frac{z_1^{|n|}-z^{|n|}}{u(\lambda)}+
z_1^{|n|}\frac{u(\lambda)-u(\lambda_1)}{u(\lambda)u(\lambda_1)}.
$$
Since $z, z_1 \in {\bf D}_1$, we have
$$
|z^{|n|}-z_1^{|n|}|=|(z-z_1)\sum_{j=0}^{|n|-1}z^{|n|-1-j}z_1^{j}|\le
|(z-z_1)\sum_{j=0}^{|n|-1}1|\le |n||z-z_1|.
$$
Moreover, we have by using (\ref{S5z(lambda)})
$$
|u(\lambda)-u(\lambda_1)| \leq \left(\max_{\nu\in I(\lambda,\lambda_1)}|u'(\nu)|\right)|\lambda-\lambda_1|,
$$
$$
|z(\lambda) - z(\lambda_1)| \leq 2 \left(1 + \max_{\nu\in I(\lambda,\lambda_1)}|u'(\nu)|\right)|\lambda-\lambda_1|.
$$
The above inequalities yield
$$
|r_0(n,\lambda)-r_0(n,\lambda_1)|\le
|\lambda-\lambda_1|\left(\frac{2|n|}{|u(\lambda)|} + \max_{\nu\in I(\lambda,\lambda_1)}|u'(\nu)|\Big(2|n| + \frac{1}{|u(\lambda_1)|}\Big)\right).
$$
Interchanging $\lambda$ and $\lambda_1$, and adding the resulting inequalities, we have
\begin{equation}
|r_0(n,\lambda)-r_0(n,\lambda_1)|\le
|\lambda-\lambda_1|(1 + |n|)M(\lambda,\lambda_1).
\label{S5roM}
\end{equation}
We also have by (1)
\begin{equation}
|r_0(n,\lambda)-r_0(n,\lambda_1)|\le
 N(\lambda,\lambda_1).
\label{S5roN}
\end{equation}
By (\ref{S5roM}) and (\ref{S5roN}), we obtain (\ref{S5r0lambda-r0lambda1}).

Using (\ref{S3hatR0(z)}) and (\ref{S5r0lambda-r0lambda1}), we have for $z \in {\bf D}_1$
\begin{equation}
\begin{split}
& \|\widehat \rho_{\alpha}\widehat q_1\, (\widehat R_0(\lambda)-\widehat
R_0(\lambda_1))\,\widehat \rho_{\alpha}\widehat q_2\|_{HS}^2 \\
& = \sum_{n,m\in
{\bf Z}} |r_0(n-m,\lambda)-r_0(n-m,\lambda_1)|^2|\widehat
\rho_{\alpha}(n)\widehat q_1(n)|^2|\widehat \rho_{\alpha}(m)
\widehat q_2(m)|^2 \\
& \le
\sum_{n,m\in {\bf Z}}
|\lambda-\lambda_1|^{2\alpha}(1+|n-m|)^{2\alpha}C_{\alpha}(\lambda,\lambda_1)^{2}
|\widehat \rho_{\alpha}(n)\widehat q_1(n)|^2|\widehat \rho_{\alpha}(m)
\widehat q_2(m)|^2 \\
& \le |\lambda-\lambda_1|^{2\alpha}C_{\alpha}(\lambda,\lambda_1)^2 \sum_{n,m\in {\bf Z}}|\widehat
q_1(n)|^2| \widehat q_2(m)|^2,
\end{split}
\nonumber
\end{equation}
which proves (\ref{holder1}).
The inequality (\ref{holder3}) follows easily from (\ref{Cholder1}).
\qed

\medskip
We study the case $d = 2$.


\begin{lemma} Let $d=2$ and $\widehat q(n)=\widehat q_1(n_1)\,\widehat q_2(n_2)$, where $\widehat q_j \in \ell^2({\bf Z})$ and $n = (n_1,n_2) \in {\bf Z}^2$. Then there exists a constant $C > 0$ such that
 \begin{equation}
\|\widehat q\,\widehat R_0(\lambda)\,\widehat q \|\leq C\|\widehat q_1\|_{\ell^2({\bf Z})}^2\|\widehat q_2\|_{\ell^2({\bf Z})}^2\big|\log \big(\lambda(\lambda - 1)(\lambda - 2)\big)\big|,
\label{S5R0estimaten=2}
\end{equation}
for all $\lambda \in \Lambda_2\cap\{|\lambda| < 3\}$.
\end{lemma}

Proof. We prove the lemma by passing it on the torus. The idea consists in reducing it to the 1-dimensional case, regarding the remaining variable as a parameter. We put
\begin{equation}
q_j(x_j) = (2\pi)^{-1/2}\sum_{n_j\in{\bf Z}}\widehat q_j(n_j)\,e^{in_jx_j},
\nonumber
\end{equation}
and define the convolution operator $q_j\ast$ by
$$
\left(q_j\ast f\right)(x) = \int_0^{2\pi}q_j(x_j-y_j)f(y)dy_j, \quad
f \in L^2({\bf T}^2),
$$
where $y = (y_1,x_2)$ if $j = 1$, $y = (x_1,y_2)$ if $j=2$.
We put
\begin{equation}
\mu = \mu(\lambda,x_2) = \lambda - h(x_2), \quad
h(x_2) = \frac{1}{2}\left(1 - \cos x_2\right),
\nonumber
\end{equation}
and define the 1-dimensional operator $A_1(\mu)$ with parameter $\mu$ by
\begin{equation}
A_1(\mu) = q_1\ast\left(h(x_1) - \mu\right)^{-1}q_1\ast =
q_1\ast\left(H_0 - \lambda\right)^{-1}q_1\ast.
\nonumber
\end{equation}
Take $f, f' \in L^2({\bf T}^2)$, and let $\widehat f, \widehat f' \in \ell^2({\bf Z}^2)$ be their Fourier coefficients. We are going to estimate
\begin{equation}
\begin{split}
C_2(\lambda) & := (\widehat q_1\,\widehat q_2\,(\widehat H_0 - \lambda)^{-1}\widehat q_1\, \widehat q_2\,\widehat f,\widehat f')_{\ell^2({\bf Z}^2)} \\
&= \left(q_1\ast q_2\ast(H_0-\lambda)^{-1}q_1\ast q_2\ast f,f'\right)_{L^2({\bf T}^2)}.
\end{split}
\nonumber
\end{equation}
Letting
$$
g = q_2\ast f, \quad g' = q_2\ast f',
$$
we  have
$$
C_2(\lambda) = \int_0^{2\pi}\left(\big(A_1(\mu)g\big)(\cdot,x_2),g'(\cdot,x_2)\right)_{L^2({\bf T}^1)}dx_2.
$$
By Lemma 5.3 (2), we obtain
$$
|C_2(\lambda)| \leq \|q_1\|^2_{L^2({\bf T}^1)}\int_0^{2\pi}
\frac{\|g(\cdot,x_2)\|_{L^2({\bf T}^1)}\|g'(\cdot,x_2)\|_{L^2({\bf T}^1)}}{|\mu(\lambda,x_2)(\mu(\lambda,x_2) - 1)|^{1/2}}dx_2.
$$
Since
$$
\|\left(q_2\ast f\right)(\cdot,x_2)\|_{L^2({\bf T}^1)} \leq
\|q_2\|_{L^2({\bf T}^1)}\|f\|_{L^2({\bf T}^2)},
$$
which follows from a simple application of Cauchy-Schwarz inequality,
we  have
$$
|C_2(\lambda)| \leq \|q_1\|^2_{L^2({\bf T}^1)}\|q_2\|^2_{L^2({\bf T}^1)}\|f\|_{L^2({\bf T}^2)}\|f'\|_{L^2({\bf T}^2)}D_2(\lambda),
$$
$$
D_2(\lambda)=
\int_{0}^{2\pi} \frac{dx_2}{|\mu(\lambda,x_2)(\mu(\lambda,x_2)-1)|^{1/2}}.
$$
Then the lemma is proved if we show for $|\lambda| < 4$
\begin{equation}
D_2(\lambda) \leq 8\left(K(\lambda) + K(\lambda - 1) + K(\lambda - 2)\right),
\label{S5C0lambdaestimate}
\end{equation}
\begin{equation}
K(\lambda) = \sqrt2\left(2 + \log\frac{3\sqrt2}{|\lambda|}\right).
\nonumber
\end{equation}

To prove (\ref{S5C0lambdaestimate}), we let
\begin{equation}
J(\lambda) = \int_0^1|(s-\lambda)s(1-s)|^{-1/2}ds,
\nonumber
\end{equation}
and first derive
\begin{equation}
D_2(\lambda) \leq 2J(\lambda) + 2J(\lambda-1).
\label{S5C0z2J}
\end{equation}
In fact, by the change of variable $s = (1-\cos x_2)/2$, we have
$$
D_2(\lambda) = 2\int_0^1|(s-\lambda)(s+1-\lambda)s(1-s)|^{-1/2}ds.
\nonumber
$$
Using  the inequality
\begin{equation}
\frac{1}{|a(a-1)|^{1/2}}=\biggr|\frac{1}{a}-\frac{1}{a-1}\biggr|^{1/2}\le
\frac{1}{|a|^{1/2}}+\frac{1}{|a-1|^{1/2}}
\label{S5aa-11/2}
\end{equation}
with $a = s+1-\lambda$,
we obtain (\ref{S5C0z2J}).

In order to compute $J(\lambda)$, we put
\begin{equation}
J_0(\lambda) = \int_0^{1/2}|s(s-\lambda)|^{-1/2}ds.
\nonumber
\end{equation}
Then we have
$$
J(\lambda)\le 2\int_0^{1/2}
\frac{ds}{|s(s-\lambda)|^{\frac12}}+2\int_{1/2}^1
\frac{ds}{|(1-s)(s-\lambda)|^{1/2}}=2J_0(\lambda)+2J_0(\lambda-1).
$$
This, combined with (\ref{S5C0z2J}), implies
\begin{equation}
D_2(\lambda) \leq 4\left(J_0(\lambda) + 2J_0(\lambda-1) + J_0(\lambda-2)\right).\label{S5C0z4J0z}
\end{equation}

Let us first consider the case $|\lambda|\geq 1$. Estimating as in (\ref{S5aa-11/2}), we have for $0 \leq s \leq 1/2$
$$
\frac{1}{|s(s-\lambda)|^{1/2}} \leq \frac{1}{|\lambda|^{1/2}}
\left(\frac{1}{|s|^{1/2}} + \frac{1}{|s-\lambda|^{1/2}}\right) \leq
\frac{1}{s^{1/2}} + \frac{1}{(1-s)^{1/2}}.
$$
Hence for $|\lambda| \geq 1$
\begin{equation}
J_0(\lambda)\leq \int_0^{1/2}
\biggr(\frac{1}{s^{1/2}}+\frac{1}{(1-s)^{1/2}}\biggr) ds =2.
\label{S5J01}
\end{equation}

 Next we consider the case $|\lambda|\leq 1$ and let $\lambda=\mu+i\nu$. If
$|\nu|\geq 4|\mu|$, we put $s=|\nu|t$, $E=1/(2|\nu|)$,
$\mu_0=\mu/|\nu|$. Note that $\sqrt 2|\nu|\ge |\lambda|\ge |\nu|$ and if  $t\geq 1/2$ then
$t-\mu_0\geq t/2$, hence $|(t-\mu_0)^2+1|^{1/4} \geq (t/2)^{1/2}$. We now compute
\begin{equation}
\begin{split}
J_0(\lambda) & = \int_0^E \frac{dt}{\sqrt t|(t-\mu_0)^2+1|^{1/4}} \\
& =\int_0^{1/2}
\frac{dt}{\sqrt t|(t-\mu_0)^2+1|^{1/4}}+ \int_{1/2}^E \frac{dt}{\sqrt
t|(t-\mu_0)^2+1|^{1/4}} \\
& \leq \int_0^{1/2} \frac{dt}{\sqrt t}+ \sqrt 2\int_{1/2}^E \frac{dt}{t}=\sqrt
2+\sqrt 2\log (2E)\le 2\sqrt 2+\sqrt 2\log {\frac{\sqrt2}{|\lambda|}}.
\end{split}
\label{S5J02}
\end{equation}
If $|\nu|\leq 4|\mu|$, we let $s=|\mu|t$,
$R=1/(2|\mu|)$, $\sigma(\mu) = \mu/|\mu|$, $\nu_0=\nu/|\mu|$. We then obtain
\begin{equation}
\begin{split}
J_0(\lambda) & =\int_0^R \frac{dt}{\sqrt t|(t-\sigma(\mu))^2+\nu_0|^{1/4}} \\
&=\int_0^{1/2}
\frac{dt}{\sqrt t|(t-\sigma(\mu))^2+\nu_0^2|^{1/4}}+ \int_{1/2}^R \frac{dt}{\sqrt t
|(t-\sigma(\mu))^2+\nu_0^2|^{1/4}} \\
& \le \int_0^{1/2} \frac{dt}{\sqrt t |t-1|^{1/2}}+ \int_{1/2}^R \frac{dt}{t-1} \\
& \leq 2 + \log2R \leq 2 + \log\frac{3\sqrt2}{|\lambda|}.
\end{split}
\label{S5J03}
\end{equation}
In view of (\ref{S5J01}), (\ref{S5J02}) and (\ref{S5J03}), we have
$$
J_0(\lambda) \leq \sqrt2\left(2 + \log\frac{3\sqrt2}{|\lambda|}\right),
$$
which, together with (\ref{S5C0z4J0z}), proves (\ref{S5C0lambdaestimate}).
\qed

\medskip
Finally we consider the case $d \geq 3$.


\begin{lemma}

(1) Let $d\geq 3$ and
$$
\widehat Q_s(n_1,\cdots,n_d)=\widehat q_1(n_1)\,\widehat q_2(n_2)\,
\widehat \rho_s(n'),
$$
 where $n' = (n_3,\cdots,n_d) \in {\bf Z}^{d-2}$, $\widehat q_1,
 \widehat q_2\in \ell^2({\bf Z})$
and $\widehat \rho_s(n') = (1 + |n'|^2)^{-s/2}$ with $0 < s \leq 1$.
 Then there exists a
constant $C_s > 0$ such that the following estimate holds:
 \begin{equation}
\|\widehat Q_s\,\widehat R_0(\lambda)\,\widehat Q_s \|\leq C_s
\|\widehat q_1\|_{\ell^2({\bf Z})}^2\|\widehat q_2\|_{\ell^2({\bf Z})}^2
\label{S5R0estimaten=3}
\end{equation}
for all $\lambda \in \Lambda_d\cap\{\lambda \in {\bf C}\,; \, |\lambda| < 2d\}$.

(2) Moreover, let $\lambda,\lambda_1\in {\bf C_\pm}\cap\{\lambda \in
{\bf C}\,; \, |\lambda| < 2d\}$ and let $g>0$ be small enough and
$$
\widehat Q_{s,g}(n_1,\cdots,n_d)=\widehat \rho_g(n_1)\widehat
Q_{s+2g}(n_1,\cdots,n_d),\ \ \ \ n\in {\bf Z}^d.
$$
 Then there exists a constant $C_{s,g} > 0$ such
that the following estimate holds true:
\begin{equation}
\label{holderd} \|\widehat Q_{s,g}\, (\widehat R_0(\lambda)-\widehat
R_0(\lambda_1))\,\widehat Q_{s,g}\|\le |\lambda-\lambda_1|^{g}
C_{s,g} \|\widehat q_1\|_{\ell^2({\bf Z})}^2\|\widehat
q_2\|_{\ell^2({\bf Z})}^2,
\end{equation}
\end{lemma}

Proof. (1) The proof is similar to the one for the previous lemma.
Let
$$
A_2(\mu)=q_1\ast q_2\ast (H_0-\lambda)^{-1}q_1\ast q_2\ast = q_1\ast q_2\ast (h_1+h_{2}-\mu)^{-1}q_1\ast q_2\ast,
$$
$$
 h_j= h(x_j) = \frac{1}{2}(1 - \cos x_j),\quad   \mu=\mu(\lambda,x') = \lambda-\sum_{j=3}^dh_j,
$$
where $x' = (x_3,\cdots,x_d)$.
For $f, f' \in L^2({\bf T}^d)$,
we put $g_s = \rho_s\ast f$, $g'_s = \rho_s\ast f'$. Then we have
\begin{equation}
\begin{split}
C(\lambda)  :=& \, (\widehat q_1\,\widehat q_2\, (\widehat H_0-\lambda)^{-1}\widehat q_1\, \widehat q_2 \,\widehat \rho\widehat f,\widehat \rho\,\widehat f') \\ = & \, (q_1\ast q_2\ast (H_0-\lambda)^{-1}q_1\ast q_2\ast\rho_s\ast f,\rho_s\ast f') \\
 =& \, \int_{{\bf T}^{d-2}}(A_2(\mu)g_s(\cdot,x'),g'_s(\cdot,x'))_{L^2({\bf T}^2)}dx'.
\end{split}
\nonumber
\end{equation}
Lemma 5.4 then implies
\begin{equation}
\begin{split}
|C(\lambda)|
& \leq C\|\widehat q_1\|_{\ell^2({\bf Z})}^2\|\widehat q_2\|_{\ell^2({\bf Z})}^2 \\
&\ \ \ \ \times\int_{{\bf T}^{d-2}}|\log(\mu(\mu-1)(\mu-2)| \|g_s(\cdot,x')\|_{L^2({\bf T}^2)}\|g'_s(\cdot,x')\|_{L^2({\bf T}^2)}dx',
\end{split}
\nonumber
\end{equation}
where $C$ is a constant independent of $\lambda \in \Lambda_d$. We now put
\begin{equation}
D_{s}(\lambda) =  \int_{{\bf T}^{d-2}}|\log(\mu(\lambda,x')(\mu(\lambda,x')-1)(\mu(\lambda,x')-2)|\|(\rho_s\ast f)(\cdot,x')\|^2_{L^2({\bf T}^2)}dx'.
\nonumber
\end{equation}
Lemma 5.5 will then be proved if we show the existence of a constant $C_s$ independent of $\lambda \in \Lambda_d\cap\{|\lambda| < 2d\}$ such that
\begin{equation}
D_{s}(\lambda) \leq C_s\|f\|_{L^2({\bf T}^d)}^2.
\label{S5Ddlambdaestimate2}
\end{equation}
We define the set ${\bf SP}$ by
\begin{equation}
{\bf SP} = \{(x_3,\cdots,x_d)\, ; \, x_j = 0 \ {\rm or} \ \pi, \ j = 3,\cdots,d\}.
\nonumber
\end{equation}
This is the set of singular points for $\mu(\lambda,x')$, since $\nabla_{x'}\mu(\lambda,x') = 0$ if and only if $x' \in {\bf SP}$. We label the points in ${\bf SP}$ by $p^{(1)},\cdots,p^{(N)}$, $N = 2^{d-2}$:
$$
{\bf SP} = \{p^{(1)},\cdots,p^{(N)}\}.
$$
For a sufficiently small $\epsilon > 0$, we put
$$
{\bf T}^{(j)}  = \{x' \in {\bf T}^{d-2}\, ; \, |x' - p^{(j)}| < \epsilon\}, \quad
1 \leq j \leq N,
$$
$$
{\bf T}^{(0)}(\lambda) = {\bf T}^{d-2}\setminus\left(\cup_{j=1}^{N}{\bf T}^{(j)}\right),
$$
and let
\begin{equation}
E^{(j)}_s(\lambda) =  \int_{{\bf T}^{(j)}}|\log(\mu(\lambda,x')(\mu(\lambda,x')-1)(\mu(\lambda,x')-2)|\|(\rho_s\ast f)(\cdot,x')\|^2_{L^2({\bf T}^2)}dx',
\nonumber
\end{equation}
\begin{equation}
E^{(j,k)}_s(\lambda) =  \int_{{\bf T}^{(j)}}|\log(\mu(\lambda,x')-k)|\|(\rho_s\ast f)(\cdot,x')\|^2_{L^2({\bf T}^2)}dx'.
\nonumber
\end{equation}
Then we have
\begin{equation}
D_s(\lambda) = \sum_{j=0}^{N}E^{(j)}_s(\lambda)
\leq \sum_{j=0}^N\sum_{k=0}^2E^{(j,k)}_s(\lambda).
\label{S5Ddlambdasum}
\end{equation}

\bigskip
We shall make use of the following version of Heinz' inequality. For the proof, see p. 232 of \cite{BiSolo87}.


\begin{prop}
Let $A, B$ be self-adjoint operators satisfying
$$
A \geq 1, \quad B \geq 0, \quad D(A) \subset D(B), \quad \|BA^{-1}\| \leq 1.
$$
Then for any $0 < \theta < 1$, we have
$$
D(A^{\theta}) \subset D(B^{\theta}), \quad
\|B^{\theta}A^{-\theta}\| \leq 1.
$$
\end{prop}

Let us estimate $E^{(j,k)}_s(\lambda)$.
 We take $\displaystyle{\theta = \frac{s}{n-2}}$, and
define self-adjoint operators $A$ and $B$ in $L^2({\bf T}^{d-2})$ by
\begin{equation}
\widehat{Af}(n') = C_0(1 + |n'|^2)^{(n-2)/2}\widehat f(n'),
\nonumber
\end{equation}
\begin{equation}
\left(Bf\right)(x') = \chi^{(j)}(x')\big|\log(\mu(\lambda,x')-k)\big|^{1/\theta}f(x'),
\nonumber
\end{equation}
where $\chi^{(j)}(x')$ is the characteristic function of the set ${\bf T}^{(j)}$, and $C_0$ is a constant to be determined later. We compute the Hilbert-Schmidt norm of $BA^{-1}$. Let $k(x')$ be the inverse Fourier image of $(1 + |n'|^2)^{-(n-2)}$. Then, up to a multiplicative constant, $BA^{-1}$ has the integral kernel
\begin{equation}
K(x',y') = C_0\chi^{(j)}(x')\big|\log(\mu(\lambda,x')-k)\big|^{1/\theta}k(x'-y').
\label{K(x'y')}
\end{equation}

One can show that for any $\alpha > 1$,
\begin{equation}
\sup_{\lambda \in \Lambda_d\cap\{|\lambda|<2d\}}\int_{{\bf T}_j}|\log(\mu(\lambda,x')-k)|^{\alpha}dx' < \infty.
\label{S5mualphaint}
\end{equation}
In fact, if $j = 0$, $\nabla_{x'}\mu(\lambda,x') \neq 0$ on ${\bf T}_0$.
Then we can take $\sum_{j=3}^d\cos x_j$ as a new variable to compute (\ref{S5mualphaint}). The case $j \neq 0$ is dealt with as follows. Suppose, for example, $p^{(j)} = (0,\cdots,0)$.
By the Morse Lemma, we can introduce new variables $y_j$, $3 \leq j \leq d$, around $p^{(j)}$ so that $\sum_{j=3}^d\cos x_j = (d-2) - \sum_{j=3}^dy_j^2$. One can then prove (\ref{S5mualphaint}) by an elementary computation. The other cases are treated similarly.

On the other hand, by Parseval's formula
$$
\int_{{\bf T}^{d-2}}|k(x')|^2dx' = \sum(1 + |n'|^2)^{-(n-2)} < \infty.
$$
Therefore $BA^{-1}$ is of Hilbert-Schmidt type, in particular, bounded. By choosing $C_0$ small enough, we have $\|BA^{-1}\| \leq 1$. Then by Proposition 5.6, $B^{\theta}A^{-\theta}$ is  bounded on $L^2({\bf T}^{d-2})$, which implies that
$$
\|E^{(j,k)}_s(\lambda)f\| \leq C_s\|f\|^2_{L^2({\bf T}^{d})},
$$
This proves (\ref{S5Ddlambdaestimate2}).

The proof of (2) repeats the  arguments from the proof of (1). \qed

\medskip
As a consequence of the above lemma, we show the following theorem.


\begin{theorem}
\label{Tsm} Let $d \geq 3$, and $\widehat H_{\gamma} = \widehat H_0+\gamma \widehat V$, where $\widehat V$ is a complex-valued potential such that $\widehat V(n)=O(|n|^{-s})$, $s>2$, as
$|n|\to \infty$, $\gamma$ being a complex parameter. Then there exists a constant $\delta > 0$ such that $\widehat H_{\gamma}$ has no eigenvalues when $|\gamma| < \delta$.
\end{theorem}

Proof. We put $\widehat Q(n_1,\cdots,n_d) = \widehat q_1(n_1)\, \widehat q_2(n_2)\, \widehat \rho(n')$, where
$$
\widehat q_1(n_1) = (1 + |n_1|^2)^{-(1+\epsilon)/4}, \quad
\widehat q_2(n_2) = (1 + |n_2|^2)^{-(1+\epsilon)/4},
\quad
\widehat \rho(n') = (1 + |n'|^2)^{-\epsilon/2}.
$$
If there exists $E \in {\bf R}$ and $\widehat f \in \ell^2({\bf Z}^d)$ such that $(\widehat H_0 + \gamma\widehat V - E)\widehat f = 0$, we have
$$
\widehat Q\,\widehat f = - \gamma\,\widehat Q\,(\widehat H_0 - E)^{-1}\,\widehat Q \cdot \widehat Q^{-1}\,\widehat V\,\widehat Q^{-1}\,\widehat f.
$$
Choosing $\epsilon > 0$ small enough, we have $\widehat Q^{-1}\,\widehat V\,\widehat Q^{-1} \in {\bf B}(\ell^2({\bf Z}^d))$. Using Lemma 5.5 and taking
$$
\frac{1}{\delta} = \|\widehat Q\,(\widehat H_0 - E)^{-1}\widehat Q\cdot
\widehat Q^{-1}\,\widehat V\,\widehat Q^{-1}\|,
$$
we obtain  Theorem 5.7. \qed


\section {Traces formulas}


\subsection{Fredholm determinant}
We shall discuss in an abstract framework in this section. Let $\mathcal H$ be a Hilbert space endowed with inner product $(\, , \, )$ and norm $\|\cdot\|$. Let ${\mathcal C}_1$ be
the set of all trace class operators on $\mathcal H$ equipped with the trace norm $\|\cdot\|_{{\mathcal C}_1}$. Recall that for $K \in {\mathcal C}_1$ and $z \in {\bf C}$, the following formula holds:
\begin{equation}
\det\,(I - zK) = \exp\left(- \int_0^z{\rm Tr}\, \big(K(1 - sK)^{-1}\big)ds\right)
\label{S6Detdefine}
\end{equation}
(see e.g. \cite{AsShi68}).
As is well-known if $A \in {\bf B}(\mathcal H;\mathcal H)$ and $B \in {\mathcal C}_1$, we have
\begin{equation}
{\rm Tr}\, (AB) ={\rm Tr}\, (BA),
\label{S6TrAB=TrBA}
\end{equation}
\begin{equation}
\det (I+ AB)=\det (I+BA).
\label{S6detAB=detBA}
\end{equation}

Suppose we are given an operator $H = H_0 + V$ on $\mathcal H$ satisfying the following conditions:

\bigskip
\noindent
{\bf (B-1)} $\ H_0$ is bounded self-adjoint.

\medskip
\noindent
{\bf (B-2)} $\ V$ is self-adjoint and trace class.

\bigskip
We put
$$
R_0(\lambda) = (H_0 - \lambda)^{-1}, \quad
\lambda \in \Lambda := \rho(H_0),
\label{S6Resolventandset}
$$
and define  $D(\lambda)$ by
\begin{equation}
D(\lambda)=\det (I+VR_0(\lambda)), \quad \lambda \in \Lambda.
\label{S6Fredholmdet}
\end{equation}


\begin{lemma}
(1) $\ D(\lambda)$ is analytic in $\Lambda$. Moreover
\begin{equation}
D(\lambda)=1+O(1/\lambda) \quad as \quad |\lambda|\to {\infty},
\label{S6Dlambdalarge}
\end{equation}
\begin{equation}
\log D(\lambda) = - \sum_{n=1}^{\infty}\frac{(-1)^n}{n}{\rm Tr}\,\left(VR_0(\lambda)\right)^n,
\label{S6logDlambda}
\end{equation}
where the right-hand side is absolutely convergent for $|\lambda| > r_0$, $r_0 >0$ being a sufficiently large constant.\\
\noindent
(2) The set $\{\lambda \in \Lambda\, ; \, D(\lambda) = 0\}$ is finite and coincides with $\sigma_d(H)$.
\end{lemma}

Proof. Letting $E_V(t)$ be the spectral decomposition of $V$, we define
$$
V^{1/2} = \int_{-\infty}^{\infty}{\rm sgn}\,t\,|t|^{1/2}dE_V(t), \quad
|V|^{1/2} = \int_{-\infty}^{\infty}|t|^{1/2}dE_V(t).
$$
There exists $r_0 > 0$ such that $\|R_0(\lambda)\| \leq C/|\lambda|$ for
$|\lambda| > r_0$, which implies
\begin{equation}
|{\rm Tr}\,\left(VR_0(\lambda)\right)| \leq C/|\lambda| \quad {\rm for}\quad |\lambda| > r_0.
\label{S6TrVR0decay}
\end{equation}
By (\ref{S6TrAB=TrBA}), taking $r_0$ large enough, we have for $|\lambda| > r_0$
\begin{equation}
\begin{split}
& {\rm Tr}\,\left(VR_0(\lambda)(1 + sVR_0(\lambda))^{-1}\right) \\
& =
{\rm Tr}\,\left(V^{1/2}R_0(\lambda)(1 + sVR_0(\lambda))^{-1}|V|^{1/2}\right) \\
&= \sum_{n=0}^{\infty}(-s)^n{\rm Tr}\,\Big(V^{1/2}R_0(\lambda)\overbrace{VR_0(\lambda)\cdots V R_0(\lambda)}^{n}|V|^{1/2}\Big) \\
&= \sum_{n=0}^{\infty}(-s)^n{\rm Tr}\,\left(VR_0(\lambda)\right)^{n+1}.
\end{split}
\nonumber
\end{equation}
We then have (\ref{S6Dlambdalarge}) and (\ref{S6logDlambda}) by (\ref{S6Detdefine}). For $\lambda \not\in \sigma(H_0)$, the eigenvalue problem $(H-\lambda)u = 0$ is equivalent to
$(I + (H_0-\lambda)^{-1}V)u = 0$, which has a non-trivial solution if and only if $\det(I + (H_0-\lambda)^{-1}V) = 0$. This proves (2) by (\ref{S6detAB=detBA}).  \qed


\begin{lemma}
The following identity holds:
\begin{equation}
\label{D2}
\log D(\lambda) =-\sum _{n \geq 1}\frac{F_n}{n}\left(\frac{1}{\lambda}\right)^{n},\quad
 F_n={\rm Tr}\,(H^n-H_0^n), \quad n\ge 1,
\end{equation}
where the right-hand side is unifomly convergent on $\{|\lambda|\ge r_0\}$ for
$r_0>0$ large enough. In particular,
\begin{equation}
\label{D3} F_1={\rm Tr}\, (V), \quad F_2={\rm Tr}\,(2VH_0+V^2).
\end{equation}
\end{lemma}

Proof. Let $R(\lambda) = (H - \lambda)^{-1}$. Take $r_0 > 0$ large enough. Then  for $|\lambda| > r_0$, we have by the resolvent equation
\begin{equation}
R(\lambda) - R_0(\lambda) = \sum_{n=1}^{\infty}(-1)^nR_0(\lambda)
\overbrace{VR_0(\lambda)\cdots VR_0(\lambda)}^{n}.
\label{S6Resolventeq}
\end{equation}
Let $F(\lambda) = \log D(\lambda)$. Since $\frac{d}{d\lambda}R_0(\lambda) = R_0(\lambda)^2$, we have by (\ref{S6logDlambda}) and (\ref{S6Resolventeq})
\begin{equation}
\begin{split}
- \frac{d}{d\lambda} F(\lambda) &= \sum_{n=1}^{\infty}(-1)^n
{\rm Tr}\,\Big(R_0(\lambda)\overbrace{VR_0(\lambda)\cdots VR_0(\lambda)}^n\Big) \\
&= {\rm Tr}\,\Big(R(\lambda) - R_0(\lambda)\Big).
\end{split}
\nonumber
\end{equation}
Using the equation
\begin{equation}
R(\lambda)= - \sum_{n=0}^{\infty}\lambda^{-n-1}H^n,
\nonumber
\end{equation}
 we obtain
\begin{equation}
F'(\lambda) = \sum_{n=0}^{\infty}\lambda^{-n-1}{\rm Tr}\,\big(H^n - H_0^n\big).
\nonumber
\end{equation}
In view of (\ref{S6Dlambdalarge}), we get (\ref{D2}).
\qed


\subsection{Spectral shift function and trace formula}
 Let $H=H_0+V$ satisfy (B-1), (B-2). Then there exists a function $\xi(\lambda)$ such that the following equality
\begin{equation}
{\rm Tr} \,(f(H)-f(H_0)) = \int_{\bf R}\xi(\lambda)df(\lambda)
\label{S6SSF}
\end{equation}
holds,
where $f$ is an arbitrary function from some suitable class. If $f$ is absolutely continuous,
then $df$ can be replaced by $f'(\lambda)d\lambda$. We call (\ref{S6SSF}) a {\it trace formula}, and
$\xi(\lambda)$ the {\it spectral shift function}
for the pair $H, H_0$.  A typical form of $\xi(\lambda)$ for our case of discrete Schr{\"o}dinger operator $\widehat H$ in \S 2  is drawn in Figure 1.

\begin{figure}
\tiny
\unitlength 1mm 
\linethickness{0.4pt}
\ifx\plotpoint\undefined\newsavebox{\plotpoint}\fi 
\begin{picture}(141.5,92.25)(0,0)
\put(1.25,86.25){\line(1,0){137.75}}
\put(2.25,52.75){\line(1,0){139}}
\put(140.25,50.25){\makebox(0,0)[cc]{$\lambda$}}
\put(58.5,49){\makebox(0,0)[cc]{$0$}}
\put(102.75,52.75){\line(0,1){1.25}}
\put(105,54.75){\makebox(0,0)[cc]{$d$}}
\put(61.75,4.25){\line(0,1){74.5}}
\put(59,77){\makebox(0,0)[cc]{$\xi$}}
\put(25.25,86.25){\line(0,1){1.25}}
\put(38.75,86.25){\line(0,1){1.25}}
\put(48.5,86.25){\line(0,1){1.25}}
\put(24.5,91.5){\makebox(0,0)[cc]{$\lambda_{1}$}}
\put(38.25,91.5){\makebox(0,0)[cc]{$\lambda_{2}$}}
\put(48,91.5){\makebox(0,0)[cc]{$\lambda_{3}$}}
\put(113.25,86.25){\line(0,1){1.25}}
\put(128.25,86.25){\line(0,1){1.25}}
\put(113,91.5){\makebox(0,0)[cc]{$\lambda_4$}}
\put(126.75,91.5){\makebox(0,0)[cc]{$\lambda_5$}}
\put(62,86.2){\linethickness{4pt}\line(1,0){41.25}}
\thicklines
\put(2.25,53){\line(1,0){23}}
\put(25.25,53){\line(0,-1){12}}
\put(25.25,41){\line(1,0){13.25}}
\put(38.5,41){\line(0,-1){10.2}}
\put(38.5,30.75){\line(1,0){10}}
\put(48.5,30.75){\line(0,-1){10.3}}
\put(48.5,20.5){\line(1,0){13.75}}
\qbezier(62,20.5)(66.125,20.25)(68.75,27.5)
\qbezier(68.75,27.5)(75.375,45.5)(81.25,42.75)
\qbezier(81.25,42.75)(88,38.5)(93.75,57.75)
\qbezier(93.75,57.75)(98.125,67.75)(103.25,68.05)
\put(103.25,68.05){\line(1,0){11.7}}
\put(115,68.05){\line(0,-1){7}}
\put(115,61.05){\line(1,0){13.75}}
\put(128.75,61.05){\line(0,-1){8.25}}
\put(128.75,53){\line(1,0){12.75}}
\end{picture}
\caption{\footnotesize The spectral shift function $\xi(\lambda)$
and five eigenvalues $\l_{1},\cdots, \l_{5}$ for $H$} \label{SF}
\end{figure}
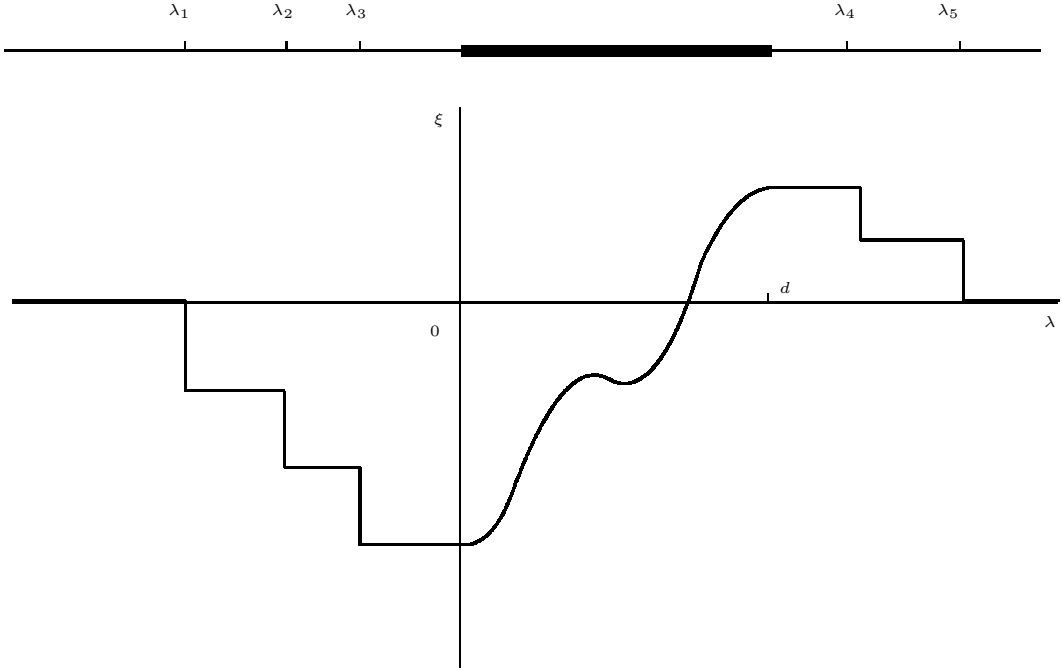

Let us recall the basic properties of $\xi(\lambda)$ (see
\cite{BY}, \cite{Ya92}).

\bigskip
\noindent
(1) {\it The following identity holds:
\begin{equation}
\log D(\lambda)=\int_{\bf R}\frac{\xi(t)}{t-\lambda}dt,\quad
\lambda\in {\bf C}_+,
\nonumber
\end{equation}
where $D(\lambda)$ is the perturbation determinant defined by (\ref{S6Fredholmdet}), and the branch of $\log D(\lambda)$ is chosen so that $\log D(\lambda)=o(1)$ as $|\lambda|\to {\infty}$,  and
$\xi(t) \in L^1({\bf R})$. We have
\begin{equation}
\xi(\lambda)=\lim_{\varepsilon\to +0}\frac{1}{\pi}\arg
D(\lambda+i\varepsilon ),\quad \ a.e. \ \lambda\in {\bf R},
\nonumber
\end{equation}
where the limit in the right-hand side exists for a.e. $\lambda\in {\bf R}$ The support of $\xi(\lambda)$ is equal to $\sigma(H)$ and
\begin{equation}
\int_{\bf R}|\xi(\lambda)|d\lambda\le \|V\|_{{\mathcal C}_1},
\nonumber
\end{equation}
\begin{equation}
\int_{\bf R}\xi(\lambda)d\lambda={\rm Tr} \,(V).
\nonumber
\end{equation}
(2) Its relation  to the S-matrix is
\begin{equation}
\det {\mathcal S}(\lambda) = e^{-2\pi i\xi(\lambda)}, \ \ for \ a.e.
\ \lambda\in \sigma_{ac}(H_0). \label{S6Smatrixxi}
\end{equation}
(3) If $H$ have $N_-\ge 0$ negative eigenvalues and
$N_+\ge 0$ positive eigenvalues, then
\begin{equation}
- N_-\le \xi(\lambda)\le N_+, \ \ for \ a.e. \ \lambda\in {\bf R}.
\nonumber
\end{equation}
(4) Suppose $H_0$ has no eigenvalues in the interval
$(a,b)\subset {\bf R}$. Assume that $\lambda_0\in (a,b)$ is an
isolated eigenvalue of finite multiplicity $d_0$ of
$H$. Then $\xi(\lambda)$ takes an integer value $n_-$ ($n_+$)  on the interval $(a,\lambda_0)$ (on the interval $(\lambda_0,b)$). Moreover, we have
\begin{equation}
\xi(\lambda_0+0)-\xi(\lambda_0-0)=-d_0.
\label{S6xiDiracdelta}
\end{equation}
(5) If $V\ge 0$, then $\xi(\lambda)\ge 0$   for all $\lambda\in {\bf
R}$.

\medskip
\noindent
(6) If $V\le 0$, then $\xi(\lambda)\le 0$  for all $\lambda\in {\bf
R}$.

\medskip
\noindent
(7) If the perturbation $V$ has  rank $N<{\infty}$, then $-N\le\xi(\lambda)\le
N$ for all $\lambda\in {\bf R}$.}

\bigskip

As will be shown in the following lemma, $F_n/n$, the Taylor coefficients of $- \log D(\lambda)$ around $\lambda = \infty$ are equal to the moments
of the spectral shift function $\xi(\lambda)$. The first two terms were computed in Lemma  6.2.
To compute the terms for $n \geq 3$, we impose the following assumption.

\bigskip
\noindent
{\bf (B-3)} {\it There exist unitary operators $S_j$ $(1 \leq j \leq d)$ such that}
\begin{equation}
\label{ConI}
 H_0= - \frac{1}{4}(S+S^*), \quad  S=\sum _{j=1}^d
S_j, \quad  S_jS_i=S_iS_j, \quad \forall\ i,j,
\end{equation}
\begin{equation}
\label{ConII}
{\rm Tr}\,\big(S_j^kV^p\big) =0,\quad  \forall \ j=1,..,d,\quad  k\neq 0,\quad
 p\geq 1.
\end{equation}

\bigskip

Note that by (\ref{ConI}) and (\ref{ConII}), we have
\begin{equation}
 {\rm Tr} \,\big(S^kV^p\big) = 0,\quad {\rm Tr}\,\big(V^p(S^{\ast})^k\big) = 0, \quad
\  k \geq 1, \quad  p \geq 1,
\label{S6SKVP0}
\end{equation}
\begin{equation}
{\rm Tr}\,\big( S^a (S^{\ast})^bV^p\big) =0, \quad  a \neq b, \quad
 p \geq 1.
 \label{S6SaSastb0}
\end{equation}


\begin{lemma}
\label{KF} Let  $H=H_0+V$ satisfy (B-1), (B-2) and (B-3).
Then
\begin{equation}
\label{KF1} F_n={\rm Tr}\,
(H^n-H_0^n)=n\int_{\bf R}\xi(\lambda)\lambda^{n-1}d\lambda,\quad
n\ge 1.
\end{equation}
In particular, letting $\tau = - 1/4$, and
\begin{equation}
 \Delta
V=\tau\sum_{i=1}^d(S_jVS_j^*+S_j^*VS_j),
\label{S6DeltaV}
\end{equation}
we have
\begin{equation}
\label{KF2} F_1={\rm Tr} \, (V),\quad F_2={\rm Tr}\,( V^2),\quad F_3={\rm
Tr}\,\big(V^3+6d\tau^2V\big),
\end{equation}
\begin{equation}
\label{KF3} F_4={\rm Tr}\, \big(V^4+8d\tau^2V^2+2\tau(\Delta V)V\big),
\end{equation}
\begin{equation}
\label{KF4} F_5={\rm Tr} \, \big(V^5+30d(2d-1)\tau^4V+10d\tau^2V^3 +5\tau(\Delta
V)V^2\big).
\end{equation}
\end{lemma}

Proof. By taking $f(\lambda) = \lambda^k$ in (\ref{S6SSF}), we get (\ref{KF1}).  In Lemma 6.2, we have proven that $F_1={\rm Tr} \, (V)$.
 $F_2$ and $F_3$ are computed by the use of (\ref{ConII}) as follows:
\begin{equation}
\begin{split}
F_2 & ={\rm Tr}\, ((H_0+V)^2-H_0^2)={\rm Tr}\, (2H_0V+V^2)={\rm Tr}\,( V^2), \\
F_3 & ={\rm Tr} \, ((H_0+V)^3-H_0^3)={\rm Tr}\, (3H_0^2V+3H_0V^2+V^3)={\rm
Tr} \, (3H_0^2V+V^3),
\end{split}
\nonumber
\end{equation}
\begin{equation}
\begin{split}
{\rm Tr}\,\big( H_0^2V\big) & =\tau^2{\rm Tr} \, \big((S^2+2SS^*+{S^*}^2)V\big)  =2\tau^2{\rm Tr}\, (SS^*V)\\
& =2\tau^2{\rm Tr} \, \big((d+\sum_{i\ne j}S_iS_j^*)V\big) =2d\tau^2{\rm Tr}\,( V).
\end{split}
\nonumber
\end{equation}

To calculate $F_4$, we first compute
\begin{equation}
\begin{split}
F_4 & ={\rm Tr}\, \big((H_0+V)^4-H_0^4)\big)\\
& ={\rm Tr}\,
(V^4+4H_0^3V+4H_0V^3+4H_0^2V^2+2H_0VH_0V).
\end{split}
\nonumber
\end{equation}
Due to (\ref{S6SKVP0}), we have
$$
{\rm Tr}\,( H_0^3V)=0,\quad {\rm Tr}\,( H_0V^3)=0,\quad
{\rm Tr} \, (H_0^2V^2) =2d\tau^2{\rm Tr} \,(V^2).
$$
Using ${\rm Tr}\,( SVSV) =0$, we get
\begin{equation}
\begin{split}
{\rm Tr} \, (H_0VH_0V) & =\tau^2{\rm Tr}\,( (S+S^*)V(S+S^*)V\big) \\
& =\tau^2{\rm Tr}\, \big((SVS^*+S^*VS)V\big)=\tau{\rm Tr} \, \big((\Delta V)V,\big),
\end{split}
\nonumber
\end{equation}
 since ${\rm
Tr}\, \big(\sum_{i\ne j}S_iVS_j^*V\big)=0$.

Finally we compute $F_5$. Firstly,
\begin{equation}
\begin{split}
F_5& ={\rm Tr}\, \big((H_0+V)^5-H_0^5\big)\\
&
={\rm Tr}\,\big(V^5+5H_0^4V+5H_0V^4+5H_0^2V^3 \\
& \ \ \ \ \ \ \ \ \ \ \ \ \ +5H_0^3V^2+5H_0^2VH_0V+5H_0VH_0V^2\big).
\end{split}
\nonumber
\end{equation}
By (\ref{S6SKVP0}), we have ${\rm Tr}\,( H_0V^4)=0,\ {\rm Tr} \, (H_0^3V^2)=0$, and ${\rm Tr}\,( H_0^2VH_0V)=0$.
We then have
\begin{equation}
\begin{split}
{\rm Tr}\,\big( H_0VH_0V^2\big) & =\tau^2{\rm Tr} \,\big((S+S^*)V(S+S^*)V^2\big) \\
& =\tau^2{\rm Tr}\,\big(
(SVS^*+S^*VS)V^2\big) \\
& =\tau{\rm Tr}\, \big((\Delta V )V^2\big)=
{\rm Tr}\,\big( H_0^4V\big)=\tau^4{\rm Tr}\,\big( (S+S^*)^4V\big) \\
& =6\tau^4{\rm Tr}\,\big(
S^2{S^*}^2V\big)=6\tau^4d(2d-1){\rm Tr}\,( V). \ \ \ \ \ \ \ \ \qed
\end{split}
\nonumber
\end{equation}

\medskip
The above lemma enables us to estimate the eigenvalues in terms of the spectral shift function.


\begin{theorem}
\label{Tes1}
 Let  $H=H_0+V$ satisfy (B-1), (B-2) and (B-3). Assume that $\sigma(H_0) = [\alpha,\beta]$, and put
\begin{equation}
 E_n = n\int_{{\bf R}\setminus[\alpha,\beta]}\xi(\lambda)\lambda^{n-1}d\lambda.
 \label{S6Gn}
\end{equation}
 Let $m_j$ be the multiplicity of $\lambda_j \in \sigma_d(H)$. Then we have for any $n \geq 0$
\begin{equation}
\label{es1-1}
\sum_{\lambda_j \in \sigma_d(H)}\, \lambda_j^n m_j=E_n.
\end{equation}
(1) If $V\ge 0$, then $\sigma_d(H) \subset (\beta,\infty)$ and
\begin{equation}
\label{es1-3} \sum_{\lambda_j \in \sigma_d(H)}\lambda_j\, m_j\leq{\rm Tr}\,( V)
\end{equation}
\begin{equation}
\label{es1-4} \sum_{\lambda_j \in \sigma_d(H)}\lambda_j^3\,m_j\leq {\rm Tr} \,\big(V^3+\frac{3d}{8}V\big).
\end{equation}
(2) If $V\le 0$, then $\sigma_d(H) \subset (-\infty,\alpha)$ and
\begin{equation}
\label{es1-6} \sum_{\lambda_j\in\sigma_d(H)}\lambda_j\,m_j\geq{\rm Tr}\,( V),
\end{equation}
\begin{equation}
\label  {es1-7} \sum_{\lambda_j\in\sigma_d(H)}\lambda_j^3\,m_j \geq {\rm Tr} \,\big(V^3+\frac{3d}{8}V\big).
\end{equation}
\end{theorem}

Proof. For small $t>0$, we define the set
$$
{\mathcal O}_t=(-{\infty},\alpha - t)\cup (\beta +t,{\infty}).
$$
If $\lambda_j \in \sigma_d(H)$, we take $\epsilon > 0$ small enough so that
the interval $I_j = (\lambda_j-\epsilon,\lambda_j+\epsilon)$ satisfies
$I_j\cap\sigma_d(H) = \{\lambda_j\}$. Then we have, by the property (\ref{S6xiDiracdelta}),
$$
\xi'(\lambda)d\lambda = - m_j\delta(\lambda - \lambda_j)d\lambda, \quad
{\rm on} \quad I_j
$$
(see Fig. 1). More precisely, see (3.28), (3.29) of \cite{BY}. We then obtain
\begin{equation}
\begin{split}
E_n(t) & := n\int_{{\mathcal O}_t}\xi(\lambda)\lambda^{n-1}d\lambda=
-\int_{{\mathcal O}_t}\xi'(\lambda)\lambda^{n}d\lambda \\
& =\sum_{\lambda_j\in {\mathcal O}_t}\lambda_j^n\,m_j \to E_n(0)=
\sum_{\lambda_j \in \sigma_d(H)}\lambda_j^n\, m_j = E_n,
\end{split}
\nonumber
\end{equation}
which proves (\ref{es1-1}).

 If $V\ge 0$, then $\lambda_j> 0$ and $\xi(\lambda) \geq 0$, which implies  $E_{2n-1}\leq F_{2n-1}$. Then (\ref{es1-1}) and (\ref{KF2}) give
(\ref{es1-3}) and (\ref{es1-4}). The proof for the case $V \leq 0$ is similar. \qed

\bigskip
\noindent
{\bf Remark 6.5} For our discrete Schr{\"o}dinger operator discussed in sections 1, 2, 3, $\widehat V = \sum_{n\in{\bf Z}^d}\widehat V_n(n)P(n)$ is trace class if
\begin{equation}
\sum_{n\in{\bf Z}^d}|\widehat V(n)| < \infty.
\nonumber
\end{equation}
The assumptions {\bf (B-1)}, {\bf (B-2)}, {\bf (B-3)} are then satisfied if we
shift our Hamiltonian $\widehat H_0$ in \S 2 by $d/2$.

\bigskip
\noindent
{\bf Remark 6.6} For the continuous model, it is well-known that $H$
has no embedded eigenvalues for the short-range perturbation. The
(non) existence of embedded eigenvalues for the discrete model is an
interesting open question.

\bigskip
\noindent
{\bf Acknowledgments.} Various parts of this paper were written during Evgeny Korotyaev's stay in the Mathematical Institute of University of Tsukuba. He is grateful to the institute for the hospitality.


\begin{thebibliography}{99}

\bibitem{AgHo76}
S. Agmon and L. H\"ormander, \textit{Asymptotic properties of solutions of differential
equations with simple characteristics}, J. d'Anal. Math., \textbf{30} (1976),
1-38.

\bibitem{AAL07}
S. Albeverio, G. F. Dell Antonio and S. N. Lakaev, \textit{The number of eigenvalues of three-particle Schr{\"o}dinger operators on lattices}, J. Phys. A : math. Theor. \textbf{40} (2007), 14819-14842.

\bibitem{ALMM06} S. Albeverio, S. N. Lakaev, K. A. Makarov, Z. I. Muminov, \textit{The threshold effects for the two-particle Hamiltonians on lattices}, Commun. Math. Phys. \textbf{262} (2006), 91-115.


\bibitem{ABG96}
W. Amrein, A. Boutet de Monvel and V. Georgescu, \textit{$C_0$-Groups, Commutator Methods and Spectral Theory of $N$-Body Hamiltonians}, Birkh{\"a}user, Progress in Math. Ser. \textbf{135}, Basel (1996)

\bibitem{AsShi68}
K. Asano and Y. Shizuta, \textit{On the theory of linear equations and Fredholm determinants}, Publ. RIMS Kyoto Univ. \textbf{3} (1968), 417-450.



\bibitem{BiSolo87}
M. S. Birman and M. Z. Solomjak, \textit{Spectral Theory of Self-Adjoint Operators in Hilbert Space}, D. Reidel Publishing Company, Dordrecht-Boston-Lancaster-Tokyo (1987).


\bibitem{BY} M. S. Birman and D. Yafaev, \textit{The spectral shift function. The papers of M. G. Krein and their further development}, St. Petersburg Math. J., \textbf{4} (1993),  883-870.



\bibitem{AMBJS99}
A. Boutet de Monvel and J. Sahbani, \textit{On the spectral properties of discrete Schr{\"o}dinger operators : (The multi-dimensional case)}, Review in Math. Phys., \textbf{11} (1999), 1061-1078.

\bibitem{B1} V. S. Buslaev, \textit{ Trace formulas for the Schr{\"o}dinger operator in a three dimensional space (Russian)}, Dokl. Akad. Nauk. SSSR, \textbf{143} (1962),1067-1070.

\bibitem{B2} V. S. Buslaev, \textit{The trace formulae and certain asymptotic
estimates of the kernel of the resolvent for the Schr{\"o}dinger
operator in three-dimensional space (Russian)}, Probl. Math. Phys. No.
I, Spectral Theory and Wave Processes, (1966) 82-101. Izdat. Leningrad
Univ. Leningrad.

\bibitem{BF}  V. S. Buslaev and L. D. Faddeev, \textit{Formulas for traces for a singular
Sturm-Liouville differential operator (Russian)}, Dokl. Akad. Nauk. SSSR, \textbf{132}
(1960), 13-16, English translation: Soviet Math. Dokl. \textbf{1} (1960),
 451-454.


\bibitem{CV}
 Y. Colin de Verdiere, \textit{Une formule de traces pour
l'op{\'e}rateur de Schr{\"o}dinger dans ${\bf R}^3$}, Ann. Sci. Ecole Norm. Sup.
(4) \textbf{14} (1981), no. 1, 27--39.



\bibitem{Fa56}
L. D. Faddeev,\textit{ Uniqueness of the inverse scattering problem}, Vestnik Leningrad Univ. \textbf{11} (1956), 126-130.

\bibitem{Fa64}
L. D. Faddeev, \textit{On the Friedrichs model in the theory of perturbations of the continuous spectrum}, Trudy Mat. Inst. Steklov \textbf{73} (1964), 292-313.

\bibitem{Fa66}
L. D. Faddeev, \textit{Increasing solutions of the Schr{\"o}dinger equations}, Sov. Phys. Dokl. \textbf{10} (1966), 1033-1035.

\bibitem{Fa76}
L. D. Faddeev, \textit{Inverse problem of quantum scattering theory}, J. Sov. Math. \textbf{5} (1976), 334-396.


\bibitem{GHS95}
F. Gesztesy, F. Holden and B. Simon, \textit{Absolute summability of the trace relation for certain Schr{\"o}dinger operators}, Commun. Math. Phys. \textbf{168} (1995), 137-161.


\bibitem{GHSZ96}
F. Gesztesy, F. Holden,  B. Simon and Z. Zhao, \textit{A trace formula for multidimensional Schr{\"o}dinger operators}, J. of Funct. Anal. \textbf{141} (1996), 449-465.


\bibitem{G}
 L. Guillop{\'e}, \textit{ Asymptotique de la phase de diffusion pour
lfop{\'e}rateur de Schr{\"o}dinger dans ${\bf R}^n$}, S{\'e}minaire EDP, 1984-1985,
Exp. No. V, Ecole Polytechnique, (1985).

\bibitem{Is03}
H. Isozaki, \textit{Inverse spectral theory}, in \textit{Topics In The Theory of Schr{\"o}dinger Operators}, eds. H. Araki, H. Ezawa, World Scientific (2003), pp. 93-143.

\bibitem{IsKu10}
H. Isozaki and Y. Kurylev, \textit{Introduction to Spectral Theory and Inverse Problems on Asymptotically Hyperbolic Manifolds}, preprint (2010).

\bibitem{JePe85}
A. Jensen and P. Perry, \textit{Commutator methods and Besov space estimates for Schr{\"o}dinger operators}, J. Operator Theory \textbf{14} (1985), 181-188.

\bibitem{Ka} N. I. Karachalios, \textit{The number of bound states for a discrete Schr{\"o}dinger operator on ${\bf Z}^N$, $N \geq 1$ lattices}, J. Phys. A 41, 455201 (2008).

\bibitem{Ka76}
T. Kato, \textit{Perturbation Theory for Linear Operators}, 2nd edition, Springer Verlag, Berlin-Heidelberg-New York (1976).

\bibitem{KhNo87}
G. M. Khenkin and R. G. Novikov, \textit{The $\overline{\partial}$-equation in the multi-dimensional inverse scattering problem}, Russian Math. Surveys \textbf{42} (1987), 109-180.

\bibitem{Mo81}
E. Mourre, \textit{Absence of singular continuous spectrum of certain self-adjoint
operators}, Commun. Math. Phys., 78 (1981), 391-408.

\bibitem{Na88}
A. Nachman, \textit{Reconstruction from boundary mesurements}, Ann. Math. \textbf{128} (1988), 531-576.

\bibitem{P}
G. Popov, \textit{ Asymptotic behaviour of the scattering phase for the
Schr{\"o}dinger operator}, C. R. Acad. Bulgare Sci. (1982) \textbf{35}(7), 885-888.

\bibitem{ReSi75} M. Reed and B. Simon, \textit{Methods of Modern Mathematical Physics, II. Fourier Analysis, Self-Adjointness}, Academic Press, New York- London, (1975).

\bibitem{R}
D. Robert, \textit{Asymptotique {\`a} grande energie de la phase de
diffusion pour un potentiel}, Asympt. Anal. \textbf{3} (1991), 301-320.

\bibitem{RoSol09} G. Rosenblum and M. Solomjak, \textit{On the spectral estimates for the Schr{\"o}dinger operator on ${\bf Z}^d$, $d \geq 3$}, Problems in Mathematical Analysis, No. 41, J. Math. Sci. N. Y. 159 (2009), No. 2, 241-263.

\bibitem{Sh} T. Shirai, \textit{A trace formula for discrete Schr{\"o}dinger operators}, Publ. RIMS Kyoto Univ. \textbf{34} (1998),  27-41.

\bibitem{SyUh87}
J. Sylvester and G. Uhlmann, \textit{A global uniqueness theorem for an inverse boundary value problem}, Ann. Math. \textbf{125} (1987), 153-169.

\bibitem{Ya92}
D. Yafaev, \textit{Mathemtical Scattering Theory: General Theory}, Amer. Math. Soc., Providence, RI, (1992).

\end{thebibliography}
\end{document}